\documentclass[a4paper,11pt]{article} 
\usepackage[utf8]{inputenc}
\usepackage[T1]{fontenc}
\usepackage{geometry}
\usepackage{amsmath,amssymb,amsfonts,amsthm,mathtools}
\usepackage{latexsym,stmaryrd,bbm}
\usepackage{nicefrac,bm}
\usepackage{graphicx}
\usepackage{tikz}
\usetikzlibrary{
  arrows.meta,
  decorations.pathreplacing,
  angles,
  quotes
}

\usepackage{xcolor}
\usepackage{hyperref}
\usepackage{algorithm}
\usepackage{algorithmicx}
\usepackage[noend]{algpseudocode}
\algrenewcommand\algorithmicrequire{\textbf{Input:}}
\algrenewcommand\algorithmicensure{\textbf{Output:}}

\usepackage{array,tabularx,booktabs}
\usepackage{microtype}
\usepackage{cleveref}
\usepackage{enumitem} 
   
\usepackage{pifont}    

\newcommand{\magsquare}{\textcolor{magenta!40}{\ding{110}}}

 \setlist[itemize,1]{label=\magsquare}  

\newtheorem{theorem}{Theorem}[section]
\newtheorem{lemma}[theorem]{Lemma}

\newtheorem{asspt}[theorem]{Assumption}
\theoremstyle{definition}

\theoremstyle{remark}

\numberwithin{equation}{section}
\numberwithin{figure}{section}
\numberwithin{table}{section}
\usepackage{subcaption}

\newcommand{\E}{\mathbb{E}}

\newcommand{\dist}{\mathbf{d}}

\newcommand{\argmin}{{\mathrm{argmin}}}

\renewcommand{\hat}{\widehat}

\usetikzlibrary{calc}
\graphicspath{{./Figures/}}

\newcommand{\email}[1]{\href{mailto:#1}{\nolinkurl{#1}}}

\title{Filtering with Randomised Observations: Sequential Learning of Relevant Subspace Properties and Accuracy Analysis}

\author{
  N.~Abedini\thanks{VU Amsterdam, Department of Mathematics, De~Boelelaan~1111, 1081~HV Amsterdam, The~Netherlands,
    \email{n.abedini@vu.nl}}
  \and
  J.~de~Wiljes\thanks{Institute for Mathematics, Ilmenau University of Technology,
    Weimarer~Str.~25, 98693 Ilmenau, Germany; and
    School of Engineering Sciences, LUT University, Finland,
    \email{jana.de-wiljes@tu-ilmenau.de}}
  \and
  S.~Dubinkina\thanks{VU Amsterdam, Department of Mathematics, De~Boelelaan~1111, 1081~HV Amsterdam, The~Netherlands,
    \email{s.b.dubinkina@vu.nl}}
}
\date{}

\begin{document}
\maketitle

\begin{abstract} 
State estimation that combines observational data with mathematical models is central to many applications and is commonly addressed through filtering methods, such as ensemble Kalman filters. In this article, we examine the signal-tracking performance of a continuous ensemble Kalman filtering under fixed, randomised, and adaptively varying partial observations. Rigorous bounds are established for the expected signal-tracking error relative to the randomness of the observation operator. In addition, we propose a sequential learning scheme that adaptively determines the dimension of a state subspace sufficient to ensure bounded filtering error, by balancing observation complexity with estimation accuracy. Beyond error control, the adaptive scheme provides a systematic approach to identifying the appropriate size of the filter-relevant subspace of the underlying dynamics. 
\end{abstract}

\section{Introduction}
\begin{figure}[ht]
    \centering

        \includegraphics[scale = 0.35]{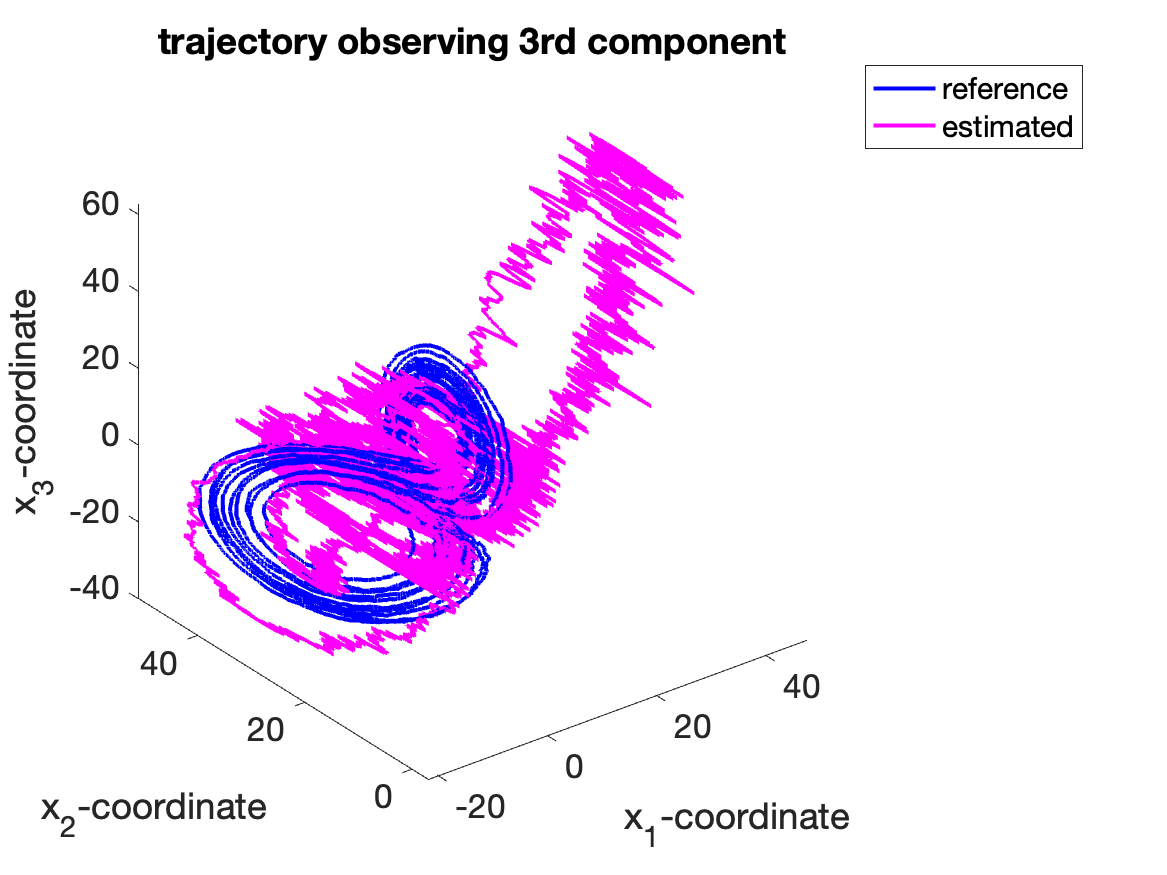}
            \includegraphics[scale = 0.35]{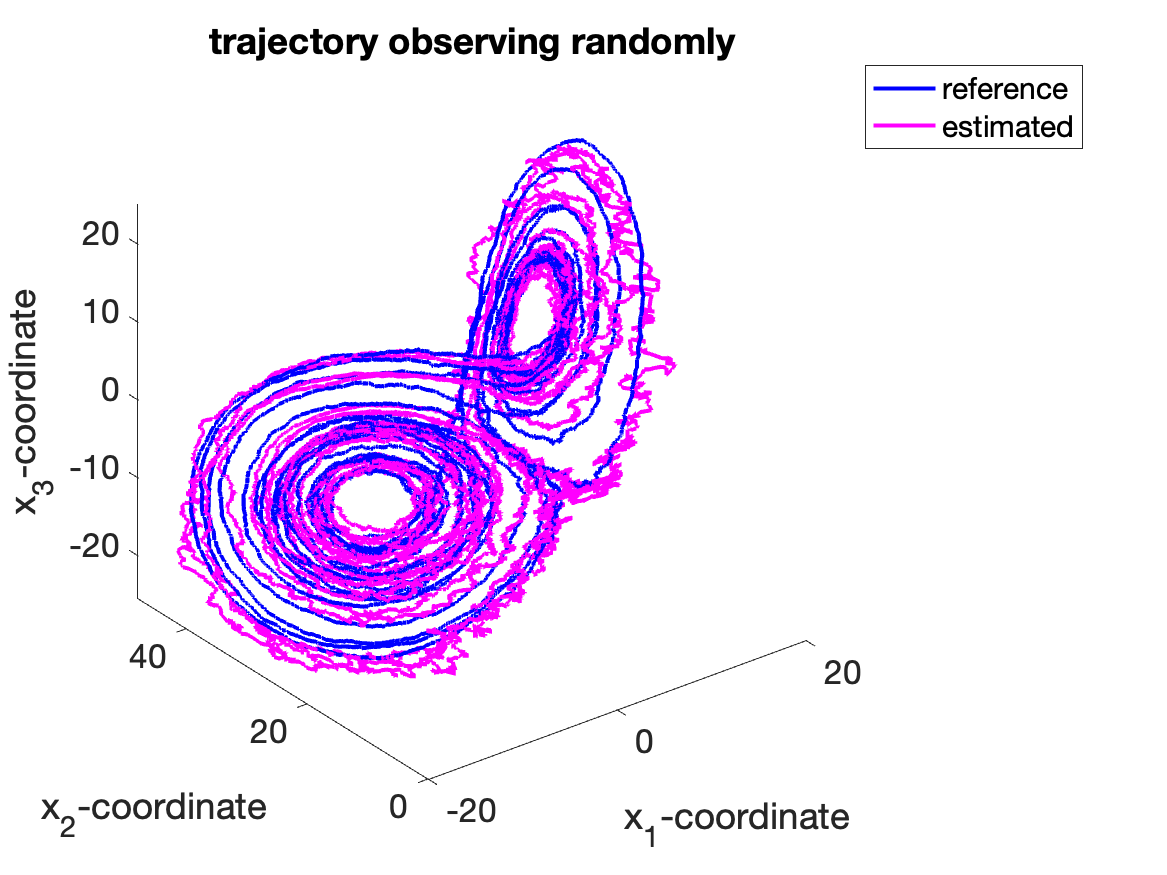}
        \caption{Solution of the Lorenz 63 system: in blue is the reference solution and in magenta is a data-assimilation solution using single noisy observation. Left: observing only the third component. Right: observing one component at random at every integration time step.}
 \label{fig:trajectory_obs_fixed}
\end{figure}

Estimating system parameters and states by combining model forecasts with observational data is a central challenge across many applications and is typically addressed with state-of-the-art filtering techniques. Two prominent application examples are climate science and numerical weather prediction \cite{lorenc1986,compo2011}, where the underlying models are typically chaotic. In such systems, even small deviations in initial conditions can lead to vastly different outcomes after only a short time. The Lorenz-63 ODE (see Section~\ref{app:L6396}) is a classical toy model that illustrates this sensitivity. Filtering accuracy depends critically on the choice of observation operator. This is illustrated in the left panel of Figure~\ref{fig:trajectory_obs_fixed}, where an Ensemble Kalman–Bucy filter fails to track a reference solution of Lorenz-63 when only the third component is observed: the filter not only diverges rapidly but even falls off the attractor. In~\cite{Carrassi-2007}, it has been shown that in order to obtain an improved estimate, the observation operator must cover the nonstable subspaces corresponding to Lyapunov vectors associated with non-negative Lyapunov exponents. In~\cite{LAW20161}, the authors numerically demonstrated, using the Lorenz-96 system, that adaptive observation operators aligned with the leading error–covariance directions can achieve estimation accuracies on the order of the observational error, even when their dimension is smaller than that of the full nonstable subspace. In this article, we analytically derive conditions on the observation operator under partial observations that ensure accurate filtering. By adaptively choosing such an operator, we guarantee a uniform bound on the estimation error. Furthermore, we prove that even if these conditions are temporarily violated, the filter can still recover a stable and relatively accurate estimate. Specifically, we consider a randomised setup in which the observed component switches according to an underlying Poisson process. This approach is illustrated numerically in the right panel of Figure \ref{fig:trajectory_obs_fixed}, where each component of the L63 model is observed randomly with probability $1/3$, subject to noise. Note that the randomisation enables accurate tracking even though the third component, which on its own is insufficient for stability, is observed with the same probability as the other two components.

Last but not least, we learn the "sufficient" size of the observation space through a novel sequential learning approach (see Section \ref{sec:sequential_learning}), which leverages a classic algorithm from stochastic multi-armed bandits (MAB) \cite{lattimore2020}. The bandit model is essentially a finite set of probability distributions, each distribution associated with an action. At each time step $t$, a learner chooses an action and receives feedback in the form of a reward, drawn randomly from the corresponding distribution. The goal in MAB is to design strategies for sequential sampling that identify the action with the highest expected reward as efficiently as possible, ideally within a relatively short time horizon.

For our theoretical considerations, we examine a deterministic Ensemble Kalman–Bucy filter (EnKBF) due to its popularity and available theoretical results on which our analysis builds. Well-posedness of discrete and continuous versions of Ensemble Kalman Filters (EnKF) in the linear setting was shown in~\cite{KellyEtAl14}, with established mean–square accuracy estimates uniform in the ensemble size. Stability of EnKF-type methods under sparse but regularly spaced observations was demonstrated in~\cite{majda08}. Subsequently, \cite{KMT15} quantified catastrophic error growth for finite ensembles when the ensemble covariance floor collapses, motivating the use of variance inflation. 
Long-time accuracy in the nonlinear setting was first obtained for the continuous-time EnKBF with complete observations in \cite{deWiljesReichStannat2018}, 
identifying a positive eigenvalue floor as the key to exponential stability. 
A Localised EnKBF was analysed for high dimensions in \cite{deWiljesTong2020}, showing how localisation can aid the accuracy estimation and improving earlier uniform bounds. Our analytical results are based on the work \cite{deWiljesTong2020} due to a natural link between localization and the definition of the observation operator. 
We do not explicitly require covariance inflation nor an observability condition on all unstable modes, properties that led to bounded long-time accuracy for a square-root EnKF in the mean-field limit, shown in \cite{SanzAlonsoWaniorek2024}. Neither do we rely on dissipativity, a property that underpins the uniform-in-time error bounds for approximate Gaussian filters applied to the 2D Navier–Stokes equations in \cite{BiswasBranicki2024}.

\subsection{Main contributions}

Here, we address the signal-tracking capabilities of a stochastic finite-ensemble localised ensemble Kalman-Bucy filter ($l$-EnKBF) with (randomised) partial observations and introduce a, within this context, novel sequential learning approach to extract relevant information about the underlying forward evolution model.

\begin{itemize}
    \item In Section \ref{sec:fixed_J_bound}, we derive error bounds for $l$-EnKBF under fixed partial observations and discuss the additional assumptions required when extending previously established results for the fully observed case.
    \item In Section \ref{sec:random_obs}, we introduce the concept of randomised partial observations in the context of filtering and extend the accuracy analysis to this more general setting, identifying the conditions under which such bounds can be derived.
    \item In Section \ref{sec:sequential_learning}, we propose a sequential learning algorithm via a multi-armed bandit approach that adaptively adjusts the number of observed components through a Poisson-randomised switching scheme. This allows for uniformly drawn observations at each switching point in a discrete-ensemble approximately Gaussian filtering framework.
    \item Finally, we validate the proposed algorithm in controlled toy settings and show that analytic subspace sizes previously determined in the literature can be approximated using the proposed algorithms (see Section~\ref{sec:sequential_learning}).
\end{itemize}

\section{Mathematical background}
We consider a continuous diffusion process $(X_t)_{t \geq 0}$ and a continuous observed signal $(Y_t)_{t \geq 0}$ governed by the following ODEs,
\begin{align} \label{eq_Diffus_full}
    \mathrm{d}X_t &= f(X_t) \mathrm{d}t + \sqrt{2}D^{1/2}\mathrm{d}W_t, \\
\label{eq_Signal_full}    \mathrm{d}Y_t &= h(t,X_t) \mathrm{d}t + R^{1/2} \mathrm{d}B_t,
\end{align}
where $dW_t$ and $dB_t$ are independent Wiener processes with diffusion covariance $D$ which, without loss of generality, will be set to the identity from here on (see \cite{deWiljesTong2020} for justification) and observation covariance $R$, respectively, having full rank. We assume that $(X_t)_{t \geq 0}$ lives on a set $\mathbb{R}^{N_x}$ and  $(Y_t)_{t \geq 0}$ lives on a set $\mathbb{R}^{N_y}$ with $N_y\le N_x$.  
Given the diffusion process $X_t$ and its associated observations $Y_{1:t}$, we consider the continuous-time filtering problem of determining the conditional distribution $\pi(X_t|\mathcal{F}(Y_{1:t}))$ 
which evolves according to the Kushner-Stratonovich equation
\[
{\rm d} \pi_t[g] = \pi_t [ f \cdot \nabla g]{\rm d}t+  \pi_t [ \nabla \cdot D\nabla g]{\rm d}t
 +  \left(\pi_t[g] - \pi_t[g]\pi_t[h]\right)^{\rm T} R^{-1}\left( {\rm d}Y_t - \pi_t[h]{\rm d}t\right),
\]
where
\begin{equation*} 
\pi_t[g] := \int_{\mathbb{R}^{N_x}} g(x) \pi_t(x){ d} x,
\end{equation*}
see~\cite{jazwinski1970} for example. Numerical methods are used to obtain approximate of this so called filtering distribution, such as sequential Monte Carlo methods, also known as particle methods~\cite{chopin2020}. Although statistically consistent, they are restricted to low-dimensional problems. On the other hand, methods such as ensemble Kalman filters, which are not consistent with respect to delivering the true posterior~\cite{Tippett2003,WhitakerHamill2002}, are notoriously known for their stability and accuracy, as shown in many high-dimensional applications, see e.g.~\cite{evensen2022,kalnay2003}. Moreover, it has been analytically shown in~\cite{deWiljesTong2020} that continuous ensemble Kalman filters, called ensemble Kalman-Bucy filters \cite{kalmanbucy1961}, are robust and accurate for abundant accurate observations. 

Here, we consider accurate \emph{partial} observations with linear time-dependent observation operator, $h(t,X_t):=H_tX_t$. We
focus on a particular family of deterministic ensemble Kalman-Bucy filters (EnKBF) where the model noise has been replaced by an interacting particle term that reflects the ensemble spread. Furthermore, we consider a localised variant which is given by the solution of a coupled system of SDEs.
The replacement of model noise is akin to the common practice of inflation, where the current spread of the ensemble is increased by a multiple of itself to maintain a certain size. This is necessary because the data assimilation update and contraction of the system often cause the particles to become very similar. Inflation is utilized to prevent filter collapse by ensuring diversity among the particles. Furthermore, since these systems are often not fully resolved, for example with respect to certain model parameters or exhibit natural variability, it is entirely reasonable to treat the associated evolution equations as a stochastic process. However, because the model noise is often poorly understood, many practitioners resort to deterministic evolution equations and compensate for the lack of variability or uncertainty in the model description through inflation. 

An ensemble $\{X^{(i)}\}_{i=1}^M$ is initially independently drawn from the distribution $\pi_0$ and then propagated forward in time by the following ODE,
\begin{equation}
\label{eq:EnKBF}
dX^{(i)}_t=f(X^{(i)}_t)dt+ P^\dagger_t(X_t^{(i)}-\overline{X}_t)dt- \frac{1}{2}P^L_t H_t^T R^{-1}(H_tX^{(i)}_tdt+H_t \overline{X}_t dt-2dY_t),\ \forall i = 1,\ldots,M,
\end{equation}
where
\[
\overline{X}_t=\frac1M \sum_{i=1}^M X^{(i)}_t,\quad P_t= \frac1{M-1} \sum_{i=1}^{M} (X^{(i)}_t-\overline{X}_t)(X^{(i)}_t-\overline{X}_t)^T,\mbox{ and }[P^L_t]_{i,j}:=[P_t]_{i,j} \phi_{i,j},
\]
with $\phi_{i,j}$ being a predefined bounded function with compact support, called localization function. 
While generally localisation is employed because spatial correlation decreases with distance in many applications, neighbouring observations are considered essential. This approach is implemented via a localisation matrix $\phi$ with 
\begin{equation}
\label{eqn:phi}
\phi_{i,j}=\rho\Big(\frac{\dist(i,j)}{r_\text{loc}}\Big),\quad\forall i,j\footnote{To simplify the notation, when
the index set is $\{1,\dots,N_x\}$ it is not specified.}, 
\end{equation}
with $r_\text{loc}$ being the so-called localisation radius and $\dist(\cdot,\cdot)$ being a distance measure between components $i$ and $j$.  Note that the localisation function $\rho$ in \eqref{eqn:phi} has to be chosen such that the diagonal entries are one and decay with distance to zero, see \cite{Gaspari1999,Stanley2021} for examples.

We assume to have a diverse ensemble in order not to consider the ensemble size effect on the error propagation. Although it is an interesting topic, it remains outside the scope of this article. Since the sample covariance matrix $P_t$ is not full rank when the number of samples is smaller than the dimension of the state space, a common scenario in classical data assimilation applications, its inverse is not defined. To address this, we replace the inverse with a diagonal approximation of the inverse i.e, the inverse of the sample covariance is replaced by the diagonal inverse \[
[P_t^\dagger]_{i,i}=[P_t]_{i,i}^{-1},\quad [P_t^\dagger]_{i,j}=0,\quad \forall i,j,\ i\neq j.
\]

\section{Bounds for signal tracking ability for fixed time}\label{sec:fixed_J_bound}
We consider the error between the ensemble mean $\overline{X}_t$ and the reference solution $X_t$,   
\begin{equation}\label{eq:trackingerror}
e_t=X_t-\overline{X}_t
\end{equation}
and aim to derive bounds for the expected error $\mathbb{E}\|e_t\|$. The comparison is particularly valid when the true posterior has a dominant mode or is approximately Gaussian with respect to a probability metric, such as the total variation distance. This assumption is often reasonable in high-dimensional data settings with relatively low observation noise, even if the underlying prior is multimodal. This brings us to the first and central assumption.
\begin{asspt}\label{asspt:error_R}
Let $\Omega = \varepsilon R^{-1}$ be a diagonal matrix, bounded by the constants 
\begin{equation}
\omega_{\min} I \preceq \Omega \preceq \omega_{\max} I,
\end{equation}
where $\varepsilon > 0$ is assumed to be small. The factor $\varepsilon$ serves as a reference order to bound the precision of the tracking ability.
\end{asspt}
In other words, we assume that the observations are reliable and that their variability is spatially independent. The latter assumption is often valid too, as in 
practice subsampling or preprocessing techniques are applied to observations to mitigate plausible dependencies. In order to be able to bound the error, we first have to ensure that the filter ensemble does not collapse or blow-up. These phenomena can be encapsulated by examining the empirical covariance, which evolves over time according to
\begin{equation}
\label{eq:evolutionP_t}
\frac{d}{dt} P_t= (F_t+F_t^T)+(P_t^\dagger P_t+P_t P_t^\dagger ) -\frac1{2\varepsilon} \left( P^{L}_tH_t^{T} \Omega H_t P_t + P_t H_t^{T} \Omega H_t P^L_t \right)
\end{equation}
with $F_t:=\frac1{M-1} \sum (X^{(i)}_t-\overline{X}_t)(f(X^{(i)}_t)-\overline f_t)^T$. The goal is to determine upper and lower bounds for $\|P_t\|_{\text{max}}$. This can be done by utilizing the fact that, for any positive semidefinite matrix, such as an empirical covariance matrix, the maximum norm is given by
\[
\|P\|_{\max} = \max_{k} [P]_{k,k},
\]
where $[P]_{k,k}$ represents the diagonal entries of the matrix $P$, corresponding to the variances in each component. We denote $\|P\|_{\min}:= \min_{k}{[P]_{k,k}}$ for its symmetry with the norm $\|P\|_{\max}$, even though it is not a norm. 
Since the first two terms on the right hand side of Equation~(\ref{eq:evolutionP_t}) are not affected by the choice of the observation operator $H_t$, they can be bounded as previously derived in~\cite{deWiljesTong2020}. 
However, for fixed $t$ the third term is highly dependent on the incoming observations.In the following, we discuss the necessary adaptations to existing bounds and proofs that account for partial observations of the signal of interest.

\subsection{Fixed observations}\label{sec:Partial_fixed_obs}
First, we focus on the component-wise error $[e_t]_k$ from Equation~\ref{eq:trackingerror} for a fixed but arbitrary state component $k$, as the accuracy is highly dependent on the pairing of observed components and the components to be updated.
We assume to have fixed in time and space, noisy yet direct\footnote{If the state space is not directly observable, it can be reformulated as a vector in the observation space, which can then be mapped back to the desired state space.} observations of state components $j\in J$, where $J\subset \{1,\dots,N_x\}$ with $N_J:=|J|$. This means that at a fixed time $t$ the observation operator $H_t$ is replaced in Equations~\ref{eq:EnKBF} and \ref{eq:evolutionP_t} by $H_J\in\mathbb{R}^{N_J\times N_x}$ with matrix entries defined by 
\begin{equation}\label{def:H_J}
[H_J]_{k,l}=\begin{cases}
1, & k=l \text{ and }k\in J;\\
0, &\text{otherwise.}
\end{cases}
\end{equation}
To ultimately bound the empirical covariance evolving according to Equation~\eqref{eq:evolutionP_t}, we need to find an upper bound for
\begin{align}
\frac{1}{2\varepsilon} \left( P^{L}_tH_J^{T} \Omega H_J P_t + P_t H_J^{T} \Omega H_J P^L_t \right).
\end{align}
Note that $\Omega$ for the particular assumptions made here is in $\mathbb{R}^{N_J\times N_J}$. Furthermore, $H^{T}_J \Omega H_J\in \mathbb{R}^{N_x\times N_x}$ and only has entries different from zero on the diagonal entries $[H^{T}_J \Omega H_J]_{j,j}$ for $j\in J$. Considering component $k$ we investigate the diagonal entries of the empirical covariance matrix:
\begin{align}
\label{eq:thirdterm}
[P_t^{L}H_J^{T} \Omega H_J P_t]_{k,k}=\sum_{i=1}^{N_x}[P_t^{L}]_{k,i}[H_J^{T} \Omega H_J]_{i,i}[P_t]_{i,k}&=\sum_{i=1}^{N_x} [P_t]_{k,i} \phi_{k,i} [H_J^{T} \Omega H_J]_{i,i} [P_t]_{i,k}\notag\\
&= \sum_{j\in J} [P_t]_{k,j} \phi_{k,j}\Omega_{j,j} [P_t]_{j,k}. \notag
\end{align}
In cases where $k$ is directly observed, the bound for the individual component can be derived in the same manner as in the fully observed settings discussed in \cite{deWiljesTong2020}. However, when $k \notin J$, noisy observations of components $j \in J$ can still be advantageous, provided that the correlation between at least one spatial components $j^{\star}$ and $k$ is sufficiently high. This brings us to the following assumption.
\begin{asspt}[Observation localisation coverage]\label{asspt:coverage}
For a localisation radius $r_\text{loc}>0$ there exists a strictly positive constant $\phi_{\min}>0$ such that for every state index \(k\in\{1,\dots,N_x\}\) one can find at least one observed index \(j^{\star}(k)\in J\) (possibly \(j^{\star}=k\)) satisfying
\[
   \phi_{k,j^{\star}}\ge\phi_{\min}.
\]
Equivalently, for 
every state index $k\in\{1,\dots,N_x\}$ there exists a nonempty localization coverage set 
defined as $\mathcal{V} = \{\forall j^{\star}\in J\mbox{ s.t. } \dist(k,j^{\star})\le r_\text{loc}\}$ for which the following holds
\[
     \max_{j\in \mathcal{V}}\phi_{k,j}
     \ge\phi_{\min},
   \quad\forall k.
\]
Furthermore let there be strictly positive constants $c_{k,j^{\star}}:=[P_t]_{k,j^{\star}}^{2}/[P_t]_{k,k}^{2}$ and  
\begin{equation}\label{eq:delta_k}
q_k:=\min_{j\in\mathcal{V}}\{c_{k,j}  \phi_{k,j}\},
\end{equation}
with  $q_k=\Theta(\kappa_{\varepsilon})$\footnote{Notation $q=\Theta(\kappa)$ is  for an equality $c_1\kappa\le q\le c_2\kappa$ with some strictly positive constants $c_1$ and $c_2$.} 
for every $k$ and each $t$.
\end{asspt}

To control the first term on the right-hand side of Equation~\eqref{eq:evolutionP_t}, it is essential to understand the role and behavior of $F_t$. This term reflects how the ensemble spread interacts with the system’s nonlinear dynamics. On the attractor, the spatial decay of interactions in the drift, characterized by the sequence $ \mathcal{F}_{\dist(i,j)}$, is typically well-behaved. This is formulated as the following assumption.
\begin{asspt}[Short-range interactions]
\label{aspt:short}
There exists a sequence of non-negative Lipschitz coefficients \( \{\mathcal{F}_k\}_{k \in \mathbb{N}} \) such that for any two state vectors \( X = [x_1, \ldots, x_{N_x}] \) and \( X' = [x'_1, \ldots, x'_{N_x}] \), the component-wise drift function \( f_i \) satisfies the bound
\[
|f_i(X) - f_i(X')| \leq \sum_{j=1}^{N_x} \mathcal{F}_{\dist(i,j)} |x_j - x'_j|, \quad \forall i.
\]
To quantify the maximal cumulative sensitivity of the system, we define the global Lipschitz constant
\begin{equation}
\label{eq:Cf}
C_\mathcal{F} := \max_{i} \sum_{j=1}^{N_x} \mathcal{F}_{\dist(i,j)}.
\end{equation}
\end{asspt}
However, if the ensemble is pushed far from the attractor, this interaction term can induce significant growth in directions not sufficiently constrained by observations. Therefore, the assumption must be interpreted as implicitly requiring that the ensemble remains sufficiently close to the attractor during the data assimilation step, where the spatial decay of the drift interaction is preserved. 
This brings us to our next crucial assumption, which says that the spatial coupling in the drift term decays sufficiently with distance and can be effectively bounded by the localization function $\phi$.

\begin{asspt}[Local Lipschitz drift aligned with localisation]
\label{asspt:drift_localization}
If the ensemble trajectories satisfy  
      \(\|X_t^{(i)}\|\le b\) for some $b>0$ for all \(i,t\) inside that ball the spatial interactions in the drift term, denoted by \( \mathcal{F}_{\dist(i,j)} \), are bounded by the localization function \( \phi \) as follows:
\[
\mathcal{F}_{\dist(i,j)} \leq C_{\mathcal{F}}  \phi_{i,j} \quad \forall i, j,
\]
where \( C_{\mathcal{F}} > 0 \) is a constant independent of \( i, j \).
\end{asspt}
Under the assumptions made so far we obtain:
\begin{align*}
\sum_{j \in J} [P_t]_{k,j} \phi_{k,j} \Omega_{j,j} [P_t]_{j,k}
&\ge [P_t]^2_{k,j^{\star}} \phi_{k,j^{\star}} \Omega_{j^{\star},j^{\star}} \quad (\text{here $j^{\star}=\argmin_{j\in\mathcal{V}}\{c_{k,j}\phi_{k,j}\}$}) \\
&\ge c_{k,j^{\star}} \phi_{k,j^{\star}} \Omega_{j^{\star},j^{\star}} [P_t]^2_{k,k}\ge c_{k,j^{\star}}  \phi_{k,j^{\star}} \omega_{\min}[P_t]^2_{k,k} \\
&\ge\min_{j\in\mathcal{V}}\{c_{k,j}  \phi_{k,j}\}\omega_{\min}[P_t]^2_{k,k}=q_k\omega_{\min}[P_t]^2_{k,k}.
\end{align*}
Using the upper bounds for the first two terms on the right-hand side of Equation~\eqref{eq:evolutionP_t} previously derived in \cite{deWiljesTong2020}, we obtain
\begin{align}\label{boundedPtineq}
    \frac{d}{dt} [P_t]_{k,k} \le -\frac{  q_k\omega_{\min}}{\varepsilon}  [P_t]^2_{k,k} +2 C_{\mathcal{F}} [P_t]_{k,k} + 2 ,
\end{align}
The following lemma provides the new upper and lower bounds of the covariance matrix, essentially indicating that, for a fixed time, the covariance should neither fully collapse nor should the spread increase beyond the specified bounds.

\begin{lemma}[Covariance bound]\label{lem:bound_Pt}
Let Assumptions \ref{asspt:error_R}, \ref{asspt:coverage}, \ref{aspt:short}, and \ref{asspt:drift_localization} hold and let the covariance matrix \( P_t \) evolve in time according to Equation~\eqref{eq:evolutionP_t}.
\begin{itemize}
    \item Then \( P_t \) satisfies the following upper bound:
\[
\|P_t\|_{\max} \leq \lambda_{\max} := \left[\left(\frac{2\varepsilon}{\omega_{\min}q_{\min}}\right)^2  C_{\mathcal{F}}^2 + \frac{8\varepsilon}{\omega_{\min}q_{\min}}\right]^{1/2},
 \quad \forall t > t_{\star} :=  \dfrac{\varepsilon}{\omega_{\min}\lambda_{\max}q_{\max}},
\]
where $q_{\min} := \min_{k} q_k$ and $q_{\max} := \max_{k} q_k$. Moreover, for all \( t > 0 \), we have 
\[
\|P_t\|_{\max} \leq \max\left\{ \|P_0\|_{\max}, \lambda_{\max} \right\}.
\]
\item
For a sufficiently small \(\varepsilon>0\), \( P_t \) satisfies the following lower bound:
\end{itemize}
\[
\|P_t\|_{\min}\ge\lambda_{\min}: = \min\Bigl\{
      \dfrac{1}{4C_{\mathcal F}^2\lambda_{\max}},
      \dfrac{\varepsilon}{C_\phi^{\star}\omega_{\max}\lambda_{\max}}
      \Bigr\}, \quad \forall t > t_{\star}+  \lambda_{\min},
 \]
where $C_\phi^{\star}:=\max_{k}\sum_{j\in J}\phi_{k,j}$. Moreover, for all $t>0$, we have  
 \[
 \|P_t\|_{\min}\ge\min\{\|P_0\|_{\min},\lambda_{\min}\}>0.
 \]
\end{lemma}
Proof of the lemma is provided in Appendix \ref{app:Bound_covariance}. We note that the constant \(q_{\min}\) influences the lower floor only via
\(\lambda_{\max}\). The error bounds obtained in~\cite{deWiljesTong2020} remain valid under the current set of assumptions and the following additional one that helps to construction of component-wise Lyapunov weights used to decouple the component errors.

\begin{asspt}
\label{aspt:diagd} 
The localization matrix $\phi$ is row-wise diagonally dominant. That is, there exists a constant $q <1$ such that for every row index $i$, the sum of the off-diagonal entries satisfies
\[
C^{(J)}_{\phi}:=\sum_{j\in J,\ j \ne i} \phi_{i,j} \leq q.
\]
\end{asspt}
\begin{lemma}\label{lem:comp_error}
 Let   
Assumptions \ref{asspt:error_R}, \ref{asspt:coverage}, \ref{aspt:short}, \ref{asspt:drift_localization}, and \ref{aspt:diagd}
hold; and let $q_{\min}$, $q_{\max}$, $\lambda_{\min}$, and $\lambda_{\max}$  denote the constants from Lemma~\ref{lem:bound_Pt}.
Let the localization weighed correlation $q_{\min}$ from Assumption \ref{asspt:coverage} be $q_{\min} = \Theta(1)$. 
Consider the error $e_t$ defined in Equation~\eqref{eq:trackingerror}, then for every fixed \(t_{0}>0\) there exists constants $c>0$ and \(C=C\bigl(q_{\min},q_{\max},\lambda_{\min},\lambda_{\max},\omega_{\min},\omega_{\max},\phi_{\min},N_J,C_{\mathcal{F}},C^{(J)}_{\phi}\bigr)>0\) such that for sufficiently small $\varepsilon$,
\begin{enumerate}
\item For all \(t> t_{0}\) and any index \(i\),
      \[
        \E[e_{t}]_{i}^{2}
           \le
            C\sqrt{\varepsilon}.
      \]

\item For any \(0<\lambda<c\varepsilon^{-1/2}\) and index \(i\),
      \[
        \sup_{t> t_{0}}
        \E\exp\bigl(\lambda [e_{t}]_i^{2}\bigr)
          \le
          2\exp\bigl(4C\lambda\sqrt{\varepsilon}\bigr).
      \]
\item For any \(T>t_{0}\) and index \(i\),
      \[
        \E_{t_{0}}\Bigl[
          \sup_{t_{0}\le t\le T}[e_{t}]^2_i
        \Bigr]
        \le
        \max_{i}[e_{t_{0}}]_i^{2}
        +C\sqrt{\varepsilon}
         \log\Bigl(\tfrac{T}{\sqrt{\varepsilon}}\Bigr).
      \]
Here $\E_{t_{0}}$ denotes conditional expectation with respect to information available at time $t_0$.
\item For any \(T>t_{0}\) and index \(i\),
      \[
        \E_{t_{0}}\Bigl[
          \max_{i}
          \sup_{t_{0}\le t\le T}[e_{t}]_i^{2}
        \Bigr]
        \le
    \max_{i}[e_{t_{0}}]_i^{2}
        +C\sqrt{\varepsilon}
         \log\Bigl(\tfrac{N_{x}T}{\sqrt{\varepsilon}}\Bigr).
      \]
\end{enumerate}
\end{lemma}

In Appendix~\ref{app:app_accuracy_fixed_J} the necessary adaptations from the original proof in \cite{deWiljesTong2020} are discussed to show that this lemma holds. In Table~\ref{tab:regimes-updated}, we provide how the covariance bound and the mean-square error scale with the order of correlation, $q_{\min}$, in terms of the observation error, $\varepsilon$. In case of high correlation, $q_{\min} = \Theta(1)$, the bounds are well behaving as in~\cite{deWiljesTong2020}, while in case of low correlation, $q_{\min} = \Theta(\varepsilon)$, though the sample covariance $P_t$ remains uniformly bounded and well-conditioned, the filter can still fail due to the mean-squared error of components scaling up as $1/\varepsilon$. However, we note that this is a component error, and thus one can still avoid filter divergence by having informative, $q_{\min}=\Theta(1)$, observations frequently. This approach is discussed in the next section.

\begin{table}[htbp]
  \centering
  \renewcommand{\arraystretch}{1.15}
  \begin{tabular}{@{}llll@{}}
\hline
    $q_{\min}$ & $\lambda_{\min}$ & $\lambda_{\max}$ &
    Component MSE\\
\hline
    $\Theta(1)$
      & $\Theta(\sqrt{\varepsilon})$
      & $\Theta(\sqrt{\varepsilon})$
      & $O(\sqrt{\varepsilon})$ \\[0.45em]
   $\Theta(\varepsilon)$
      & $\Theta(\varepsilon)$
      & $\Theta(1)$
      & $O(\varepsilon^{-1})$ \\[0.45em]
\hline
  \end{tabular}
    \caption{Asymptotic spread and error scales for two representative regimes
           $q_{\min}= \Theta(\varepsilon^{\nu})$ with $\nu = 0$ and $\nu =1$.  The balanced case
           ($\nu=0$) matches the original analysis.}
  \label{tab:regimes-updated}
\end{table}

\subsection{Randomised observations}\label{sec:random_obs}
\begin{figure}
    \centering
\begin{center}
\begin{tikzpicture}[scale=0.4]
    \draw[thick, black, rounded corners=1.5pt, ->] (-12.5,6.5) to[out=0,in=180] (-11,6.5); 
    \draw[thick, black, rounded corners=1.5pt, ->] (-12.5,5.5) to[out=0,in=180] (-11,5.5); 
    \draw[thick, black, rounded corners=1.5pt, ->] (-12.5,4.5) to[out=0,in=180] (-11,4.5); 
    \draw[thick, black, rounded corners=1.5pt, ->] (-12.5,3.5) to[out=0,in=180] (-11,3.5); 
    \draw[thick, black, rounded corners=1.5pt, ->] (-12.5,2.5) to[out=0,in=180] (-11,2.5); 
    \draw[thick, black, rounded corners=1.5pt, ->] (-12.5,1.5) to[out=0,in=180] (-11,1.5); 
    \draw[thick, black, rounded corners=1.5pt, ->] (-12.5,0.5) to[out=0,in=180] (-11,0.5); 

    \draw[thick, black, rounded corners=1.5pt, ->] (-10.5,6.5) to[out=0,in=180] (-3.5,1.5); 
    \draw[thick, black, rounded corners=1.5pt, ->] (-10.5,5.5) to[out=0,in=180] (-3.5,0.5); 
    \draw[thick, black, rounded corners=1.5pt, ->] (-10.5,4.5) to[out=0,in=180] (-3.5,6.5); 
    \draw[thick, black, rounded corners=1.5pt, ->] (-10.5,3.5) to[out=0,in=180] (-3.5,3.5); 
    \draw[thick, black, rounded corners=1.5pt, ->] (-10.5,2.5) to[out=0,in=180] (-3.5,5.5); 
    \draw[thick, black, rounded corners=1.5pt, ->] (-10.5,1.5) to[out=0,in=180] (-3.5,4.5); 
    \draw[thick, black, rounded corners=1.5pt, ->] (-10.5,0.5) to[out=0,in=180] (-3.5,2.5); 

    \draw[thick, black, rounded corners=1.5pt, ->] (-3.5,6.5) to[out=0,in=180] (4.5,0.5); 
    \draw[thick, black, rounded corners=1.5pt, ->] (-3.5,5.5) to[out=0,in=180] (4.5,2.5); 
    \draw[thick, black, rounded corners=1.5pt, ->] (-3.5,4.5) to[out=0,in=180] (4.5,6.5); 
    \draw[thick, black, rounded corners=1.5pt, ->] (-3.5,3.5) to[out=0,in=180] (4.5,5.5); 
    \draw[thick, black, rounded corners=1.5pt, ->] (-3.5,2.5) to[out=0,in=180] (4.5,4.5); 
    \draw[thick, black, rounded corners=1.5pt, ->] (-3.5,1.5) to[out=0,in=180] (4.5,3.5); 
    \draw[thick, black, rounded corners=1.5pt, ->] (-3.5,0.5) to[out=0,in=180] (4.5,1.5); 

    \draw[thick, black, rounded corners=1.5pt, ->] (4.5,6.5) to[out=0,in=180] (11.5,4.5); 
    \draw[thick, black, rounded corners=1.5pt, ->] (4.5,5.5) to[out=0,in=180] (11.5,1.5); 
    \draw[thick, black, rounded corners=1.5pt, ->] (4.5,4.5) to[out=0,in=180] (11.5,6.5); 
    \draw[thick, black, rounded corners=1.5pt, ->] (4.5,3.5) to[out=0,in=180] (11.5,0.5); 
    \draw[thick, black, rounded corners=1.5pt, ->] (4.5,2.5) to[out=0,in=180] (11.5,5.5); 
    \draw[thick, black, rounded corners=1.5pt, ->] (4.5,1.5) to[out=0,in=180] (11.5,2.5); 
    \draw[thick, black, rounded corners=1.5pt, ->] (4.5,0.5) to[out=0,in=180] (11.5,3.5); 

    \draw[thick, black, rounded corners=1.5pt, ->] (12,6.5) to[out=0,in=180] (13.5,6.5); 
    \draw[thick, black, rounded corners=1.5pt, ->] (12,5.5) to[out=0,in=180] (13.5,5.5); 
    \draw[thick, black, rounded corners=1.5pt, ->] (12,4.5) to[out=0,in=180] (13.5,4.5); 
    \draw[thick, black, rounded corners=1.5pt, ->] (12,3.5) to[out=0,in=180] (13.5,3.5); 
    \draw[thick, black, rounded corners=1.5pt, ->] (12,2.5) to[out=0,in=180] (13.5,2.5); 
    \draw[thick, black, rounded corners=1.5pt, ->] (12,1.5) to[out=0,in=180] (13.5,1.5); 
    \draw[thick, black, rounded corners=1.5pt, ->] (12,0.5) to[out=0,in=180] (13.5,0.5); 

    \filldraw[color=magenta!20] (-11,6) -- (-10,6) -- (-10,7) -- (-11,7) -- cycle; 
    \filldraw[color=magenta!20] (-11,5) -- (-10,5) -- (-10,6) -- (-11,6) -- cycle; 
    \filldraw[color=gray!20] (-11,4) -- (-10,4) -- (-10,5) -- (-11,5) -- cycle; 
    \filldraw[color=gray!20] (-11,3) -- (-10,3) -- (-10,4) -- (-11,4) -- cycle; 
    \filldraw[color=magenta!20] (-11,2) -- (-10,2) -- (-10,3) -- (-11,3) -- cycle; 
    \filldraw[color=magenta!20] (-11,1) -- (-10,1) -- (-10,2) -- (-11,2) -- cycle; 
    \filldraw[color=magenta!20] (-11,0) -- (-10,0) -- (-10,1) -- (-11,1) -- cycle; 

    \draw[color=gray, very thick] (-11,0) -- (-10,0) -- (-10,7) -- (-11,7) -- cycle;
    \foreach \y in {1,...,6} {
        \draw[color=gray, very thick] (-11,\y) -- (-10,\y);
    }
    
    \node[color=black] at (-10.5,-1.5) {\footnotesize{$x^{(i)}_{t}$}};
    
    \filldraw[color=magenta!20] (-4,6) -- (-3,6) -- (-3,7) -- (-4,7) -- cycle; 
    \filldraw[color=gray!20] (-4,5) -- (-3,5) -- (-3,6) -- (-4,6) -- cycle; 
    \filldraw[color=magenta!20] (-4,4) -- (-3,4) -- (-3,5) -- (-4,5) -- cycle; 
    \filldraw[color=magenta!20] (-4,3) -- (-3,3) -- (-3,4) -- (-4,4) -- cycle; 
    \filldraw[color=magenta!20] (-4,2) -- (-3,2) -- (-3,3) -- (-4,3) -- cycle; 
    \filldraw[color=gray!20] (-4,1) -- (-3,1) -- (-3,2) -- (-4,2) -- cycle; 
    \filldraw[color=magenta!20] (-4,0) -- (-3,0) -- (-3,1) -- (-4,1) -- cycle; 
    
    \draw[color=gray, very thick] (-4,0) -- (-3,0) -- (-3,7) -- (-4,7) -- cycle;
    \foreach \y in {1,...,6} {
        \draw[color=gray, very thick] (-4,\y) -- (-3,\y);
    }

    \node[color=black] at (-3.7,-1.5) {\footnotesize{$x^{(i)}_{t+\tau_1}$}};
    
    \filldraw[color=magenta!20] (4,6) -- (5,6) -- (5,7) -- (4,7) -- cycle; 
    \filldraw[color=gray!20] (4,5) -- (5,5) -- (5,6) -- (4,6) -- cycle; 
    \filldraw[color=gray!20] (4,4) -- (5,4) -- (5,5) -- (4,5) -- cycle; 
    \filldraw[color=magenta!20] (4,3) -- (5,3) -- (5,4) -- (4,4) -- cycle; 
    \filldraw[color=magenta!20] (4,2) -- (5,2) -- (5,3) -- (4,3) -- cycle; 
    \filldraw[color=magenta!20] (4,1) -- (5,1) -- (5,2) -- (4,2) -- cycle; 
    \filldraw[color=magenta!20] (4,0) -- (5,0) -- (5,1) -- (4,1) -- cycle; 

    \draw[color=gray, very thick] (4,0) -- (5,0) -- (5,7) -- (4,7) -- cycle;
    \foreach \y in {1,...,6} {
        \draw[color=gray, very thick] (4,\y) -- (5,\y);
    }
    
    \node[color=black] at (4.8,-1.5) {\footnotesize{$x^{(i)}_{t+\tau_2}$}};
    
    \filldraw[color=magenta!20] (11,6) -- (12,6) -- (12,7) -- (11,7) -- cycle; 
    \filldraw[color=magenta!20] (11,5) -- (12,5) -- (12,6) -- (11,6) -- cycle; 
    \filldraw[color=magenta!20] (11,4) -- (12,4) -- (12,5) -- (11,5) -- cycle; 
    \filldraw[color=magenta!20] (11,3) -- (12,3) -- (12,4) -- (11,4) -- cycle; 
    \filldraw[color=magenta!20] (11,2) -- (12,2) -- (12,3) -- (11,3) -- cycle; 
    \filldraw[color=gray!20] (11,1) -- (12,1) -- (12,2) -- (11,2) -- cycle; 
    \filldraw[color=gray!20] (11,0) -- (12,0) -- (12,1) -- (11,1) -- cycle; 

    \draw[color=gray, very thick] (11,0) -- (12,0) -- (12,7) -- (11,7) -- cycle;
    \foreach \y in {1,...,6} {
        \draw[color=gray, very thick] (11,\y) -- (12,\y);
    }
    
    \node[color=black] at (11.8,-1.5) {\footnotesize{$x^{(i)}_{t+\tau_3}$}};

    \draw[color=violet, very thick] (-11,4) -- (-10,4) -- (-10,5) -- (-11,5)--(-11,4); 
    \draw[color=violet, very thick] (-4,4) -- (-3,4) -- (-3,5) -- (-4,5)--(-4,4); 
    \draw[color=violet, very thick] (4,4) -- (5,4) -- (5,5) -- (4,5)--(4,4); 
    \draw[color=violet, very thick] (11,4) -- (12,4) -- (12,5) -- (11,5)--(11,4); 

    \begin{scope}[shift={(-11, -3)}] 
        \filldraw[color=magenta!20] (-0, -1) rectangle (1, -2);
        \node[color=black, anchor=west] at (1.2, -1.5) {\footnotesize{observed components $x^{(i)}_{s}(j)$ with $j\in J$}};
        
        \filldraw[color=white] (0, 0) rectangle (1, 1);
        \draw[color=violet, very thick] (0,-0.4 ) rectangle (1, 0.6);
        \node[color=black, anchor=west] at (1.2, 0.1) {\footnotesize{component $x^{(i)}_{s}(r)$ that is being updated}};
             \node[] at (11.5, 10) {\footnotesize{$f(x^{(i)}_{t+\tau_1})$}};

    \end{scope}
    
\end{tikzpicture}
\end{center}
    \caption{The graphic illustrates the randomisation mechanism at the switching points $\tau_i$, which are governed by a Poisson process with rate $\lambda$. At each switch, a subset of filter components is randomly selected, and subsequent updates are performed using this fixed set until the next switching point. The figure highlights how changes in the observed states (coloured in magenta) may influence the filter update (purple frame indicates to be updated component) through this stochastic selection process.}
    \label{fig:visualisation_switches}
\end{figure}
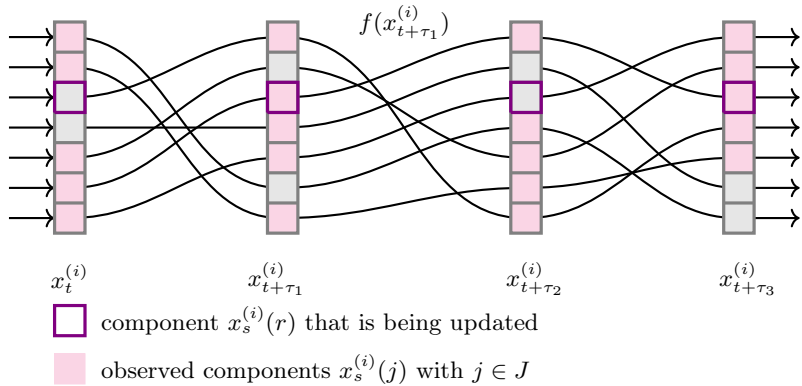

We now replace the fixed index set $J_t\subset\{1,\dots N_x\}$
by a sequence of identical independent random variables
\[
   \mathcal{J}_t = J_t \subset \{1,\dots,N_x\},
   \qquad t \ge 0,
\]
each drawn uniformly from all sets of cardinality $N_J$, i.e.\
\[
   \mathbb{P}\big(\mathcal{J}_t = J_t\big)
   =\binom{N_x}{N_J}^{-1}.
\]
with $N_{J}$ fixed (Figure \ref{fig:visualisation_switches} illustrates the effects of this randomisation within the context of filtering). For large $\binom{N_x}{N_J}$ computing all subsets and selecting one uniformly at random is impractical. Thus, efficient algorithms such as reservoir sampling~\cite{tille2006sampling,Vitter1985,Knuth1997} are typically employed to sample a subset of size $N_J$ uniformly at random without generating all subsets. An exemplary run for \(N_J=1\), with a switch occurring at every time step for the Lorenz 63 system, is visualized in the right panel of Figure \ref{fig:trajectory_obs_fixed}. Although the trajectory is recovered only noisily, it remains stable at this high switching frequency (here every time step) and with the given number of observables in the state space. However, if the switching frequency were significantly reduced (see next subsection), it could easily lead to the blow-up observed in the left panel of Figure \ref{fig:trajectory_obs_fixed}.

Next, we make the following assumption.
\begin{asspt}[Probabilistic coverage]
\label{ass:Event_c_probab}
For each time $t\ge 0$ and a component $k$ denote
\begin{equation}\label{eq:_event_E}
   \mathcal{E}_t^{k}:=
   \Bigl\{\exists j\in J_t:
          q_{k}(t) = \Theta(1)\mbox{ and }\phi_{k,j}\ge\phi_{\min}\Bigr\},
\end{equation}
then for each time $t\ge 0$
\[
   \mathbb{P}_{J\sim\mathcal J}\bigl(\mathcal{E}_t:=\cap_{k}\mathcal E_t^{k}\bigr)
   \;\ge\;1-\delta 
\]
with $\delta=O\bigl(\varepsilon^{\tfrac32+\eta}\bigr)$, $0<\varepsilon<1$ and $\eta>0$.
\end{asspt}
Since the observation operator $H_{J}$ is time-dependent, it is likely that for each component $k$ that is not observed we can find a component $j$ that is being observed sufficiently over time assuming $N_J$ is large enough relative to $N_x$.
We note that size $N_J$ influences probability $\delta$, e.g. in the extreme case $N_x = N_J$ we obtain $\delta = 0$ due to the trivial choice $j^{\star}=k$ yielding $q_{k}=1$ for every $k$. Therefore, it is natural to select $N_J$ based on the order of $\delta$. 
In addition to the theoretical result below, in Section~\ref{sec:sequential_learning} we explore a sequential approximative learning to estimate $N_J$. Since Assumption \ref{ass:Event_c_probab} replaces the deterministic
coverage Assumption~\ref{asspt:coverage} for the remainder, we derive the following result. 
\begin{lemma}[Accuracy random observation]\label{lem:random_obs_acc}
Let Assumptions~\ref{asspt:error_R}, \ref{aspt:short}, \ref{asspt:drift_localization}, \ref{aspt:diagd},
and \ref{ass:Event_c_probab} hold. 
Consider the error $e_t$ defined in Equation~\eqref{eq:trackingerror}, then for every fixed \(t_{0}>0\) there exist a constant $D>0$ such that for sufficiently small $\varepsilon$ and for all
$t>t_0$,
\begin{enumerate}
\item For any index $i$,
      \[
         \E_{J\sim\mathcal J}\Bigl[
               \E\bigl[e_i^{2}(t)\bigr]
            \Bigr]
            \le
            D_t^{N_J}\sqrt{\varepsilon}.
      \]

\item For any \(0<\eta<1\) and any index $i$,
      \[
         \mathbb{P}_{J\sim\mathcal J}\Bigl(
           \E[e_i^{2}(t)]
           \le \tfrac{D_t^{N_J}}{\eta}\sqrt{\varepsilon}
         \Bigr)
         \ge 1 - \eta  .
      \]

\end{enumerate}
\end{lemma}

The proof of the lemma can be found in Appendix~\ref{app:proofrandobs_fixedt}. We also numerically validate the error order of the continuous filter. We perform numerical experiments with the L63 model and an integration time step $d t = 0.000005$ to diminish the effect of numerical integration error. The empirical expected error is computed by averaging over 50 independent runs, where the observations are switched at every time step, and each run consists of $1{,}000{,}000$ time steps. The error order is then evaluated for different values of $\varepsilon$ for the L63 model. The numerically obtained results agree with the theoretical order and thus provide computational confirmation of the proven results.
\begin{figure}[ht]
  \centering
  \includegraphics[width=0.8\textwidth]{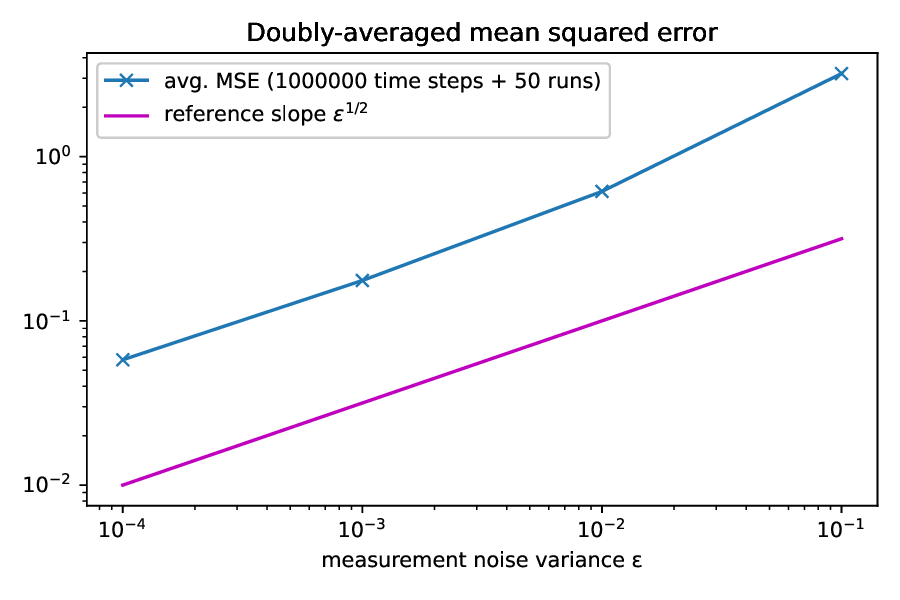}
  \caption{Numerical verification of the theoretical bound for the simulations of the continuous filter for the L63 system with $d t=0.000005$ in doubly-averaged mean squared error for different observation error order $\epsilon$.}
  \label{fig:mse_vs_eps}
\end{figure}
\section{Sequential learning of number of components \texorpdfstring{$N_J$}{NJ}}\label{sec:sequential_learning}
Since it is computationally, as well as analytically, non-trivial, particularly in high dimensions, to derive the number of meaningful components to fulfil Assumption~\ref{ass:Event_c_probab}, we propose a sequential learning strategy to learn an appropriate $N_J $. 
\subsection{Poisson–sampled \texorpdfstring{$J_k$}{Jk}}
Let $\{\tau_k\}_{k\ge1}$ be the jump times of a Poisson process $\mathcal{N}_t$ with intensity $\lambda>0$, independent of all other noise sources.  At each jump, we draw a new index set $J_k:=J_{\tau_k}\subset\{1,\dots,N_x\}$ i.i.d. from the law
$\mathcal{J}$ introduced in Section \ref{sec:random_obs} and keep it fixed until the next jump:
\[
      J(t)=J_k
      \quad\text{for }t\in[\tau_k,\tau_{k+1}).
\]
For $t\in[\tau_k,\tau_{k+1})$ each sample $X^{(i)}$ obeys the continuous-time l-EnKBF driven by the current mask $J_k$, i.e., 
\begin{equation*}\label{eq:filter_between_jumps}
  dX^{(i)}_t= f\bigl(X^{(i)}_t\bigr)dt
       +P_t^{\dagger}\bigl(X^{(i)}_t-\overline X_t\bigr)dt-\frac12 P^{L}_{t}
          H_{J_k}^{T}R^{-1}
          \bigl(H_{J_k}X^{(i)}_tdt + H_{J_k}\overline X_tdt-2dY_t
          \bigr),\ \forall i=1,\dots,M.
\end{equation*}
At jump times $t=\tau_k$, the mask $J_k$ inducing a new observation operator is instantaneously replaced by $J_{k+1}$ while everything else remains continuous. This form of update of the observation operator is more realistic for practical applications than the continuous switches we considered so far. 

\subsection{Stochastic multi-armed bandit}
In order to learn $N_J$ sequentially, we model the problem as a stochastic multi-armed bandit, where each arm corresponds to a specific choice of $N_J$. The reward is defined as a performance for that specific choice and is treated as a random variable due to the underlying components being sampled at random. Thus, their contributions can vary markedly even when $N_J$ is held constant. Moreover, the filter-update performance of a given component set can drift over time due to system dynamics and other nonlinear factors. Although these effects imply that rewards are generally not independent across arms, the component sampling itself is uniform and independent of all previous choices or previously observed set sizes. Sequential algorithms that explicitly model such dependencies do exist, but for our initial numerical study we deliberately employ a simple algorithm, namely an UCB algorithm~\cite{Auer2002,Auer2002b} that assumes temporal independence and independence among arms. The UCB algorithm is described in Algorithm~\eqref{alg:UCB1} in Appendix \ref{app:algorithms}.

Although the correlation across arms is certainly non-negligible, we often sample using a step size tailored to the state dimension $N_x$, especially in large state dimensions where one would not obtain a sufficient number of plays for each $N_J$ otherwise. Reward values for unsampled $N_J$ can then be interpolated from the sampled ones. If the step size between arms is chosen sufficiently large, the resulting correlation becomes small enough to justify assuming independence across arms. At first we need to define the corresponding reward function: 
\begin{equation}\label{def:reward_function}
\mathcal{R}_{a_{\tau}} = \underbrace{\beta \,\kappa}_\text{coverage} \;-\; \underbrace{\alpha \,\frac{N_J}{N_x}}_{\text{penalty term preventing too large $N_J$}} -  \underbrace{\gamma\,\frac{\text{trace}(P_t)}{N_x}}_{\text{spread}}
\end{equation}
with $\alpha$, $\beta$, and $\gamma$ being hyperparameters and $a_\tau=N_J$.

\begin{figure}[ht]
\begin{center}
\begin{tikzpicture}[
  transform shape,
  every node/.style={minimum size=0.8cm-\pgflinewidth, outer sep=0pt}
]
  \def\cols{12}
  \def\rows{12}
  \def\visibleRows{8}  

  \draw[step=1,thin] (0.5,0.5) grid ({\cols+0.5},{\visibleRows+0.5});

  \coordinate (C) at (4.5,4.5);
  \node[fill=magenta!40] at (C) {};

  \foreach \dx/\dy in {0/1,0/-1,1/0,-1/0,1/1,1/-1,-1/1,-1/-1}
    \node[fill=blue!20] at ($(C)+(\dx,\dy)$) {};

  \foreach \dx/\dy in {0/2,0/-2,2/0,-2/0,
                       1/2,2/1,-1/2,-2/1,
                       1/-2,2/-1,-1/-2,-2/-1,
                       2/2,2/-2,-2/2,-2/-2}
    \node[fill=blue!20] at ($(C)+(\dx,\dy)$) {};

  \foreach \x/\y in {1.5/1.5,3.5/3.5,5.5/2.5,7.5/4.5,9.5/6.5,11.5/7.5,
                     6.5/7.5,8.5/5.5,10.5/3.5,2.5/5.5,6.5/1.5}
    \draw[magenta,very thick] (\x-0.5,\y-0.5) rectangle (\x+0.5,\y+0.5);

\end{tikzpicture}

\def\legendScale{0.8}         

\begin{tikzpicture}[
  scale=\legendScale,
  transform shape,
  legendsquare/.style={
      draw,
      minimum width=0.3cm,
      minimum height=0.3cm,
      inner sep=0pt},
  every node/.style={anchor=west,font=\scriptsize}
]
  \matrix[column sep=8mm, row sep=2mm]{
    \node[fill=magenta!40,legendsquare] {}; &
    \node[xshift=2mm]{to be updated cell}; &

    \node[fill=blue!20,legendsquare] {}; &
    \node[xshift=2mm]{spatially correlated w.r.t.\ cell}; &

    \node[magenta,very thick,legendsquare] {}; &
    \node[xshift=2mm]{observed cells}; \\
  };
\end{tikzpicture}

\end{center}
    \caption{Visualization of the number of components that must be examined to estimate coverage. To update the magenta cell at the centre of the purple cells, only three observed cells need to be considered, because they lie within the localisation radius of the cell being updated.}
    \label{fig:coverage_alg_visulization}
\end{figure}
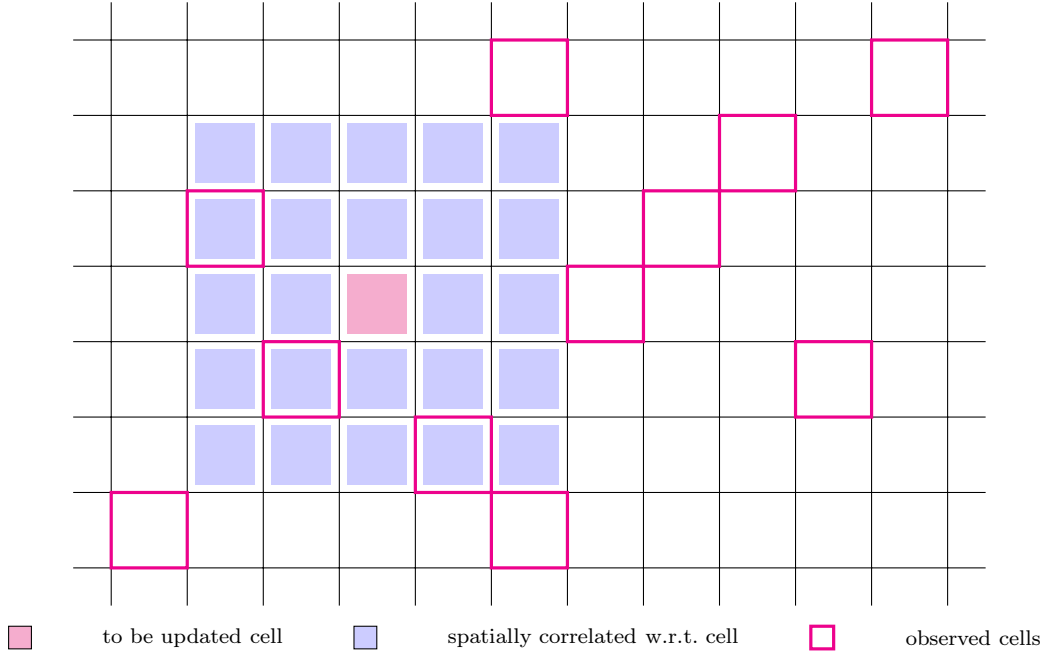

Its leading term performs a heuristic verification of the Assumption~\ref{asspt:coverage} and can be computed using Algorithm~\eqref{alg:appr_Coverage}. 
To limit computational cost, this test is carried out only within the localisation radius, because the localisation kernel vanishes beyond that range, rendering distant state components irrelevant to the variable currently being updated (see Figure \ref{fig:coverage_alg_visulization} for visualization). While the coverage ensures that the current choice of $N_J$ is sufficiently large, we introduce a penalty term to keep $N_J$ as small as possible yet still adequate (in the spirit of information criteria / Occam's razor). Note that the balance between hyperparameters $\alpha$ and $\beta$ is crucial. The reward is further enhanced by measures an estimate of the spread of the ensemble, which enters as second penalisation term in order to detect an occurring blow up. Ideally, the hyperparameter $\gamma$ is selected to preserve the spread while fine-tuning it in anticipation of the impending blow-up. Furthermore, we incrementally update the mean‑reward as follows
\begin{equation*}
\hat{\mathcal{R}}_{N_J} \;=\; 
\hat{\mathcal{R}}_{N_J} + \frac{\mathcal{R}_t - \hat{\mathcal{R}}_{N_J}}{\text{plays}_{N_J}}.
\end{equation*}

Using the proposed sequential learning algorithm we estimate a sufficient $N_J$ for the Lorenz 96 model, for the model description see Appendix~\ref{app:L6396}. Instead of a continuous l-EnKBF, we employ a discrete LETKF~\cite{WhitakerHamill2002,Tippett2003}. A rigorous connection between this discrete-time filter and its continuous-time counterpart has been demonstrated, for example in~\cite{LangeStannat2021_CMS}. Although one can write down a fully continuous formulation, any practical implementation demands a discrete-time grid anyway. Finally, the two-step update strategy (analysis and forecast) used here is numerically more robust and mirrors the schemes routinely adopted in numerical weather prediction. The proposed sequential learning algorithm with a discrete-time LETKF is given in Algorithm~\ref{alg:dacycle}. 

\begin{figure}[ht]
  \centering
  \begin{subfigure}{0.48\textwidth}
    \centering
    \includegraphics[width=\linewidth]{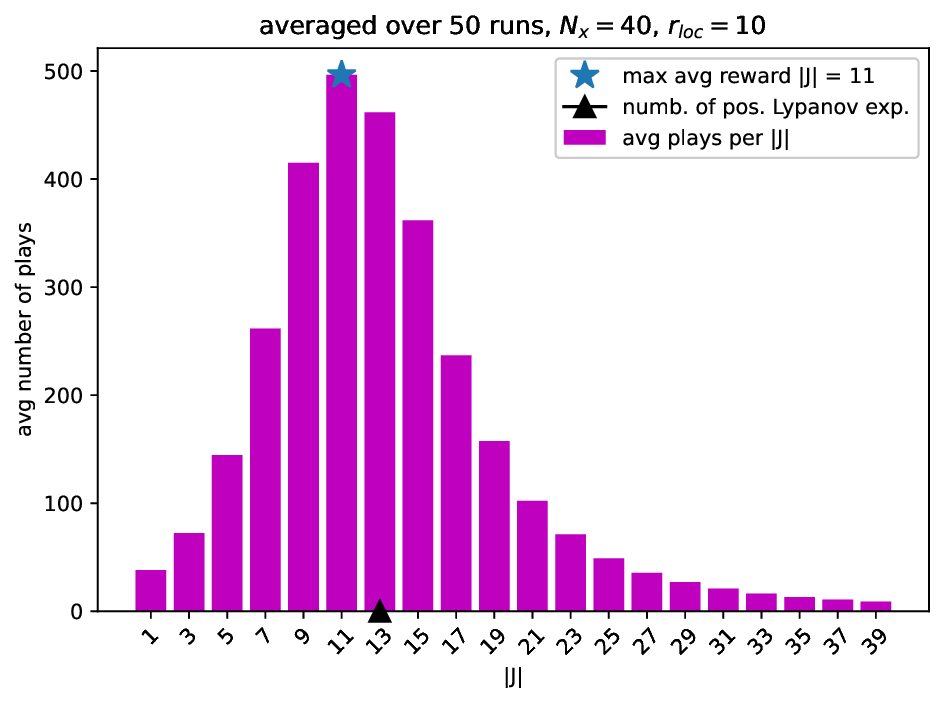}
    \caption{$N_x=40$.}
  \end{subfigure}\hfill
  \begin{subfigure}{0.48\textwidth}
    \centering
    \includegraphics[width=\linewidth]{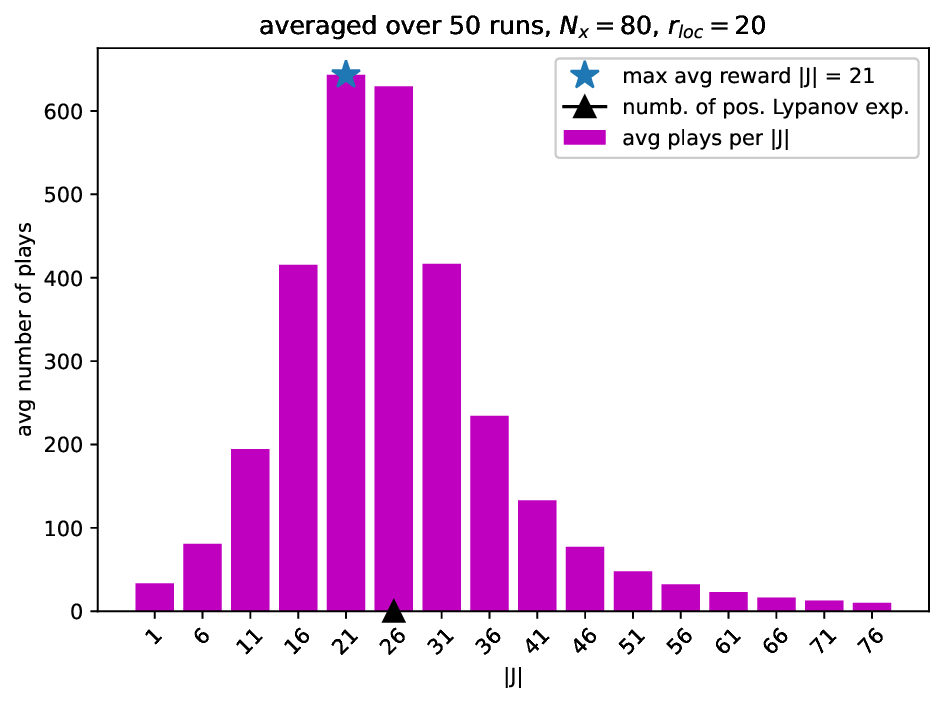}
    \caption{$N_x=80$.}
  \end{subfigure}

  \vspace{0.6em}

  \begin{subfigure}{0.48\textwidth}
    \centering
    \includegraphics[width=\linewidth]{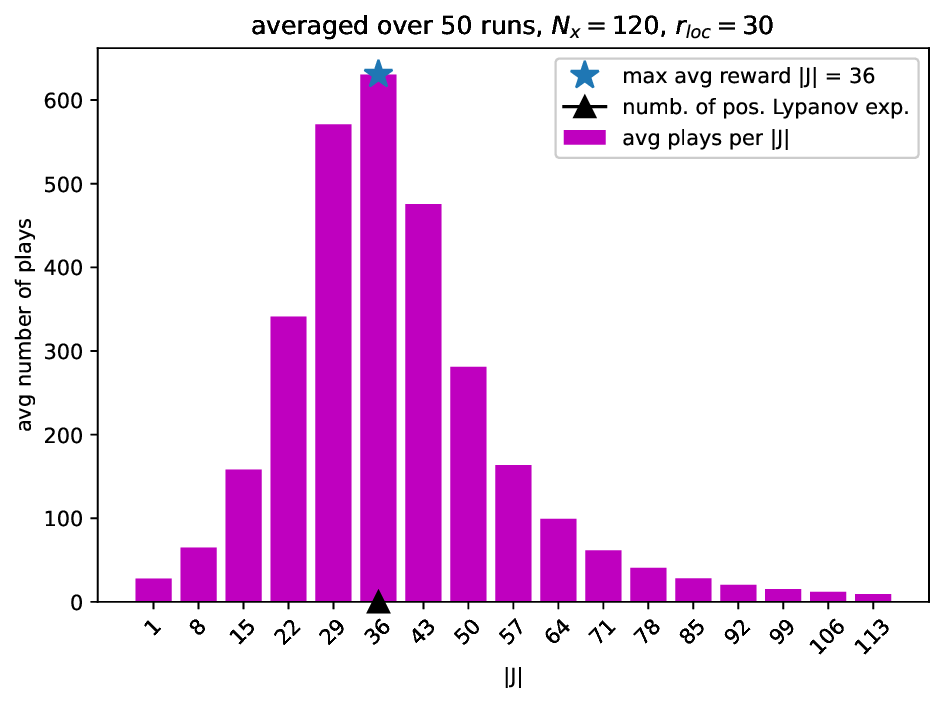}
    \caption{$N_x=120$.}
  \end{subfigure}\hfill
  \begin{subfigure}{0.48\textwidth}
    \centering
    \includegraphics[width=\linewidth]{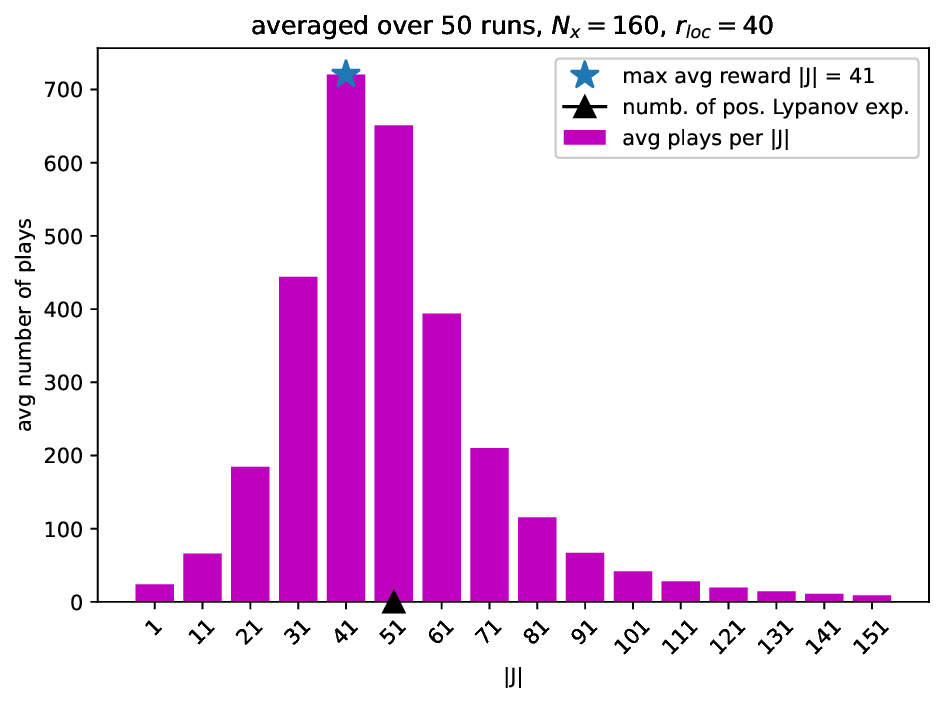}
    \caption{$N_x=160$.}
  \end{subfigure}
 \caption{Plots display number of times different values of $N_J$ being used
  in the filter for dimensions $N_x\in\{40,80,120,160\}$, the resulting optimal $N^\ast_J$, and the corresponding number
  of positive Lyapunov exponents.}\label{fig:results_J_L96}
\end{figure}

We fix the frequency $\lambda$ to $1000$ since this corresponds to a switch in every integration time step (note that the computations are robust to $\lambda$ being much smaller as even $\lambda=5$ is resulting in a very high probability to switch).
We perform numerical experiments with different state dimensions, namely $N_x\in\{40,80,120,160\}$. 
To reduce correlation across arms and computational complexity, $N_J$ belongs to a subset of $N_x$. 
That is, for $N_x=40$, we learn $N_J\in \{1:2:N_x\}$ and for $N_x=120$, we learn $N_J\in \{1:7:N_x\}$. The variables used for the computations can be found in the Table~\ref{tab:settings_runs_L96_J}. The learned $N_J$, averaged over 50 runs and across four different dimensions of the L96 system, consistently lies just below or on the known number of positive Lyapunov exponents for that dimension, as can be seen in Figure~\ref{fig:results_J_L96}. This close agreement indicates that the learned $N_J$ provides a reliable, data-driven estimate of the number of dynamically relevant components that need to be observed, consistent with analytical expectations.  
 
\section{Outlook}
The provided analysis and the newly proposed sequential learning ansatz opening the way toward a new approach to efficient data assimilation methods in complex, high-dimensional systems and is filling a gap in the rigorous analysis of a class of continuous filters. Looking ahead, several promising research directions emerge. First, the choice of reward functions in the sequential learning framework deserves closer examination, as it critically determines what is considered optimal. Moreover, identifying or even designing an algorithm that fully accounts for temporal and spatial dependencies remains an open challenge, alongside the task of establishing rigorous regret bounds that explicitly incorporate such dependencies. Furthermore, a Bayesian formulation of sequential learning appears particularly attractive, as it would naturally integrate prior information indicating which observation state sizes or specific regions of the state space are more relevant. This would allow the results to be extended beyond learning the appropriate size alone, toward also identifying areas of particular importance within the state space.
At the same time, several practical limitations must be addressed, including sensitivity to hyperparameters, the challenge of estimating regret under limited samples, and the computational complexity inherent in high-dimensional settings. In addition, it will be valuable to explore how these algorithms can inform underlying experimental design questions especially in scenarios involving extreme events that build up sequentially over time with inherent delays, or in settings where estimation must already take place online as data arrives, rather than in a batch-processing framework.
  
\section*{Acknowledgments}
The research of JdW and NA has been partially funded by the Deutsche Forschungsgemeinschaft (DFG)- Project-ID 318763901 - SFB1294. Furthermore, JdW has been supported by the European Union under the Horizon Europe Research \& Innovation Programme (Grant Agreement no. No 101188131 UrbanAIR). Additionally JdW gratefully acknowledges funding by the Carl-Zeiss foundation within the project KI-MSO-O. The completion of the manuscript, including the corresponding proofs and numerical results, was made possible through the visitors’ programme of the Amsterdam Center for Dynamics and Computation of the Vrije Universiteit Amsterdam, which JdW and SD would like to gratefully acknowledge.
\bibliographystyle{plain}
\bibliography{references_randomisedPaper}

\section{Appendix}
For improved clarity and readability, all proofs are presented in the appendix. Further technical details that are too extensive for inclusion in the main text are provided in the subsequent sections.

\subsection{Lorenz 63 and 96 models}\label{app:L6396}
The Lorenz system used to generate the trajectories shown in Figures~\ref{fig:trajectory_obs_fixed} is the classical Lorenz 63 model, a simplified representation of atmospheric convection dynamics. It is defined by the following system of three coupled nonlinear ordinary differential equations:
$$\frac{dx}{dt} = \sigma(y - x), \quad \frac{dy}{dt} = x(\rho - z) - y, \quad \frac{dz}{dt} = xy - \beta z,$$
where \( \sigma \), \( \rho \), and \( \beta \) are positive constants representing, respectively, the Prandtl number, the Rayleigh number, and a geometric factor. In our simulations, we use the standard chaotic parameter values
$\sigma = 10$, $\rho = 28$, and $\beta = 8/3$.
This system is known for its sensitive dependence on initial conditions and was first introduced in the seminal work of Lorenz~\cite{lorenz1963deterministic}. Since we incorporate the drift term in a stochastic differential equation, we perform the numerical integration with the Euler–Maruyama scheme. 

The Lorenz 96 (L96) model is a canonical high–dimensional toy atmosphere introduced by Lorenz \cite{lorenz1996predictability}.  
It consists of $N_x$ cyclically coupled scalar variables \(X_1,\dots,X_N\) that evolve according to
\begin{align}
\frac{dX_i}{dt} =\bigl(X_{i+1}-X_{i-2}\bigr)\,X_{i-1} - X_i + F,
\qquad i = 1,\dots,N_x,
\label{eq:lorenz96}
\end{align}
with periodic indexing $X_{i+N_x}\equiv X_i$.  Here we will use the standard chaotic configuration $F = 8$. The quadratic advection term
$(X_{i+1}-X_{i-2})X_{i-1}$ mimics nonlinear energy transfer, the linear damping $-X_i$ represents dissipation, and the constant forcing $F$
maintains the system out of equilibrium.

\subsection{Proof for sample covariance bound}\label{app:Bound_covariance}
Below find the proof of Lemma \ref{lem:bound_Pt} from Section \ref{sec:Partial_fixed_obs}.
\begin{proof}
First, we prove the result on the upper bound.
Our proof builds significantly on the previous work \cite{deWiljesTong2020} which is based on well established results for constant Riccati equations \cite{kucera1973review}. In particular, we employ auxiliary Lemma A.3, which asserts that if, at every time step, the maximal component of the solution to a differential equation is bounded by a quadratic function with a negative coefficient for the $x^2$ term, then it is bounded above by the difference 
\[
\Delta_{\varepsilon}:= x_1-x_2=
\frac{2\varepsilon}{\omega_{\min}q_k}\left(C_{\mathcal{F}}^2 + \frac{2\omega_{\min}q_k}{\varepsilon}\right)^{1/2}
\]
between the roots 
\[
x_{1,2} = \frac{\varepsilon}{\omega_{\min}q_k}\left[C_{\mathcal{F}} \pm \left(C_{\mathcal{F}}^2 + \frac{2\omega_{\min}q_k}{\varepsilon} \right)^{1/2} \right]
\]
of the quadratic function given in \eqref{boundedPtineq}, or by the initial condition.

To obtain the threshold time after which the maximal component of the solution to a differential equation is bounded by $\Delta_{\varepsilon}$,
we need to bound 
\[
\exp\left(\frac{\omega_{\min}q_k}{\varepsilon}t\Delta_{\varepsilon}\right)
\]
below by a constant (as shown in \cite{deWiljesTong2020}). Thus, for $t> \varepsilon/(\omega_{\min}q_k\Delta_{\varepsilon})$ we have that 
$[P_t]_{k,k}$ is bounded by $\Delta_{\varepsilon}$.
To extend these results to the maximum norm, we 
denote $q_{\min} := \min_{k} q_k$ and $q_{\max}:= \max_{k} q_k$ with $q_k$ defined in eq. \eqref{eq:delta_k} in Assumption \ref{asspt:coverage}.  
Then we have that 
\[
\|P_t\|_{\max} \leq \lambda_{\max}
:= \left[\left(\frac{2\varepsilon}{\omega_{\min}q_{\min}}\right)^2  C_{\mathcal{F}}^2 + \frac{8\varepsilon}{\omega_{\min}q_{\min}}\right]^{1/2},\quad\forall  t > t_\star:= \dfrac{\varepsilon}{\omega_{\min}\lambda_{\max}q_{\max}}.
\]
Moreover, for all \( t > 0 \),
\[
\|P_t\|_{\max} \leq \max\left\{ \|P_0\|_{\max}, \lambda_{\max} \right\}.
\]
Next, we prove the result on the lower bound. Let $k$ be such that $\|P_t\|_{\min}=[P_t]_{k,k}$.
Define
\[
   x(t):=-[P_t]_{k,k}\le 0,
   \qquad
   z_t:=\|P_t\|_{\max}>0.
\]
We need to bound each term in the component evolution equation \eqref{eq:evolutionP_t}. We focus on the third term again as the other bounds can be directly taken from the original proof in \cite{deWiljesTong2020}. Due to the fact that only indices in $J\subset\{1,\dots,N_x\}$ are observed we have
\begin{align}
\bigl[P_t^{L}H_J^{T}\Omega H_J P_t\bigr]_{k,k}
  &=
  \sum_{j\in J}\phi_{k,j}\Omega_{j,j}[P_t]_{k,j}^{2}\le
   \omega_{\max}\sum_{j\in J}[P_t]_{k,j}^{2}\phi_{k,j}
\le
   \omega_{\max}\sum_{j\in J}[P_t]_{k,k}[P_t]_{j,j}\phi_{k,j}
\nonumber\\
&\le
   \omega_{\max}[P_t]_{k,k}\|P_t\|_{\max}
   \sum_{j\in J}\phi_{k,j}
=
   \omega_{\max}C_\phi^{(J)}[P_t]_{k,k}\|P_t\|_{\max}, 
\label{eq:PO-upper_forlower}
\end{align}
where
\[
  C_\phi^{(J)}:=\sum_{j\in J}\phi_{k,j}
  \le
  C_\phi:=\sum_{i=1}^{N_x}\phi_{k,i}.
\]
Bound \eqref{eq:PO-upper_forlower} injects the mixed term, \(-xz_t\beta^{(J)}/\varepsilon\),  into the comparison function \(g(x,t) = \alpha\sqrt{-xz_t} - xz_t\beta^{(J)}/\varepsilon  - \gamma - 1\) for Lemma A.3 Claim 3 in \cite{deWiljesTong2020}, with
\[
  \beta^{(J)}
  =
C_\phi^{(J)}\omega_{\max},
\]
by substituting \(x:=-[P_t]_{k,k}\) and $z_t: = \|P_t\|_{\max} $ into 
$\dfrac{\mathrm dx}{\mathrm dt}
  \le g(x,t)$ and 
yielding 
\[
  x(t)\le -c_{\varepsilon}\mbox{ for } t>\frac{\varepsilon}{\omega_{\min}q_k\Delta_{\varepsilon}} + c_{\varepsilon},\mbox{ with }c_{\varepsilon}:=\min\Bigl\{
      \frac{\gamma^2}{9\alpha^2\Delta_{\varepsilon}},
      \dfrac{\gamma\varepsilon}{3\beta^{(J)}\Delta_{\varepsilon}}
      \Bigr\},
  \]
 $\gamma = 3$, and $\alpha = 2C_{\mathcal F}$.
 To extend this result to the minimum norm, 
 we denote $C_\phi^{\star}: = \max_{k}\sum_{j\in J}\phi_{k,j}$. Then we have that
  \[
  \|P_t\|_{\min}\ge \lambda_{\min}:=
  \min\Bigl\{
      \dfrac{1}{4C_{\mathcal F}^2\lambda_{\max}},
      \dfrac{\varepsilon}{\omega_{\max}\lambda_{\max}C_\phi^{\star}}
      \Bigr\},\quad\forall t > t_\star + \lambda_{\min}.
\]
 Moreover, for all $t>0$,
 \[
 \|P_t\|_{\min}\ge\min\{\|P_0\|_{\min},\lambda_{\min}\}>0.
 \]
Note that the early-time estimate remains unchanged.
\end{proof}

\subsection{Proof of accuracy bounds for fixed individual components fixed \texorpdfstring{$J$}{J}}\label{app:app_accuracy_fixed_J}
Below the adaptations to the original proof in are discussed that verify Lemma \ref{lem:comp_error} in Section \ref{sec:Partial_fixed_obs}.
\begin{proof}
Since the proofs of Lemma \ref{lem:comp_error} mirror those in \cite{deWiljesTong2020}, we restrict ourselves to the main differences. In the proof of Lemma B.1 in \cite{deWiljesTong2020}, which forms the basis for all subsequent bounds, the component-wise evolution must be adapted to the fact that we observe only the subset of components with indices $j\in J$, i.e.,
\[
  d[e_t]_i^{2}
  \le
     \Big(
       -\alpha_t[e_t]_i^{2}
       -2\varepsilon^{-1}
  \sum_{j\in J}[P_t\circ\tilde\phi]_{ij}[e_t]_i[e_t]_j
       +\sum_{j\in J,\ j\ne i}\mathcal{F}_{\dist(i,j)}|[e_t]_j|^{2}
       +\beta_t
     \Big)dt
   +d[\mathcal M_t]_i .
\] with $\tilde{\phi}:=\phi-I$ where $I$ is the identity matrix. Considering the analysis term, the mixed term ($j\neq i$ ) approximations remain exactly as in the fully observed case. i.e., 
\[
  -2\phi_{ij}P_{ij}[e_t]_i[e_t]_j
  \le \phi_{ij}\bigl(P_{jj}[e_t]_j^{2}+P_{ii}[e_t]_i^{2}\bigr),
\]
 The diagonal part has to be adapted to 

\[
  -2P_{ii}\phi_{ii}[e_t]_i^{2}
     +\sum_{j\in J,\ j\ne i}\phi_{ij}P_{ii}[e_t]_i^{2}
     \le (q-2)P_{ii}[e_t]_i^{2}< -P_{ii}[e_t]_i^{2}
\]
using that $q<1$. Hence the evolution of $[e_t]_i^{2}$ becomes
\[
  \mathrm d[e_t]_i^{2}
     \le
     \Bigl(
       \sum_{j\in J,\ j\ne i}
         (\mathcal{F}_{\dist(i,j)}
          +\varepsilon^{-1}\phi_{ij}P_{jj})[e_t]_j^{2}
        -\alpha_t[e_t]_i^{2}
        -\varepsilon^{-1}P_{ii}[e_t]_i^{2}
        +\beta_t
     \Bigr)\mathrm dt
     +\mathrm d[\mathcal M_t]_i .
\]
After taking the $v^{(i)}$–weighted sum of $E_t:=[[e_t]_1^2,\dots, [e_t]_{N_x}^2 ]$ for a Lyapunov function $E_t^{(i)}:=\langle v^{(i)}, E_t \rangle$, we again obtain the following inequality:
\[
   \mathrm dE_t^{(i)}
   \le
     \bigl(-\alpha_t+C_{\mathcal{F}}\bigr)E_t^{(i)}\mathrm dt
     +\beta_t\mathrm dt
     +\mathrm d\mathcal M_t^{(i)},
\]
where the weighting factors $v^{(i)}$ are defined in Lemma C.1 in~\cite{deWiljesTong2020}, which still provides
\(\sum_{j\in J\ j\neq k}\phi_{k,j}v^{(i)}_{k}\le v^{(i)}_{k}\) for all $k$. 
Then we have the following  constants
\[
  \begin{aligned}
    &\textbf{Growth coefficient:}\quad
    \alpha_t = 1-2\varepsilon^{-1}\|P\|_{\min},\\
    &\textbf{Precision bound:}\quad
        \beta_t
          =
          C_{\mathcal{F}}^{2}\|P_t\|_{\max}
          +2
          +\varepsilon^{-1}
          \omega_{\max}
           (C^{(J)}_\phi)^{2}\|P_t\|_{\max}^{2},\quad\mbox{where }C_\phi^{(J)}(i):=\sum_{j\in J}\phi_{i,j}.\\
 \end{aligned}
\]
Apart from these replacements every subsequent step (Gr\"onwall inequality,
exponential moments, pathwise maxima) is identical to the fully
observed proof. One simply propagates the modified constants
\(\omega_{\max},C_\phi^{(J)},\phi_{\min}\) through the same inequalities.

\end{proof}

\subsection{Proof for randomised observation fixed \texorpdfstring{$t$}{t} error}\label{app:proofrandobs_fixedt} 
Below that proof for Lemma \ref{lem:random_obs_acc} from Section \ref{sec:random_obs} is outlined. 

\begin{proof}
Since bounds for the event $\mathcal{E}$ have already been derived in the previous section, we now analyse the expectation, specifically, the error in terms of the two complementary events $\mathcal{E}$ as defined in equation~\eqref{eq:_event_E}, and its complement
$\left(\mathcal{E}_t^k\right)^{\mathrm{c}} = \left(\mathcal{E}_t^k\right)^{\mathrm{c}_1} \cup \left(\mathcal{E}_t^k\right)^{\mathrm{c}_2}$, where
\[
\left(\mathcal{E}_t^k\right)^{\mathrm{c}_1} 
= 
\left\{
  \forall\, j \in J_t:\; 
  q_k =\Theta(\varepsilon^\nu)\mbox{ with }0<\nu\le 1\footnote{We do not consider the case of $q_k >1$ since it does not exist in typical applications. In case it does, it remains highly improbable and thus only a slight adaption to the proof is required.} \mbox{ and } 
  \phi_{k,j} \ge \phi_{\min}
\right\}
\]
and
\[
\left(\mathcal{E}_t^k\right)^{\mathrm{c}_2} 
= 
\left\{
J_t \cap \mathcal{V}_k = \varnothing
\right\} .
\]

From Lemma~\ref{lem:comp_error}, we have
\[
\mathbb{E}\left[[e_t]_{i}^2 \mid \mathcal{E}_t\right] \le C \sqrt{\varepsilon}\quad \mbox{and}\quad\mathbb{E}\left[[e_t]_{i}^2 \mid \mathcal{E}^{\mathrm{c}_1}_t\right] \le A \varepsilon^{-\nu}.
\]
This allows us to decompose the total expected squared error:
\begin{align*}
\mathbb{E}_{J \sim \mathcal{J}}  \mathbb{E}[e_t]_{i}^2 
&= \mathbb{E}\left[[e_t]_{i}^2 \mid \mathcal{E}_t\right]  \mathbb{P}[\mathcal{E}_t] 
+ \mathbb{E}\left[[e_t]_{i}^2  \mid \mathcal{E}_t^\mathrm{c} \right] \mathbb{P}[\mathcal{E}_t^\mathrm{c}] \\
&= C_t^{N_J} \sqrt{\varepsilon} (1 - \delta) +  
A_t^{N_J} \varepsilon^{-\nu} \delta+
B^{N_J}_t \delta \\
&\le  \sqrt{\varepsilon}\left(C_t^{N_J} -C_t^{N_J} \delta + A_t^{N_J}\varepsilon^{-3/2}\delta + B_t^{N_J}\varepsilon^{-1/2}\delta \right)\\
&\le \sqrt{\varepsilon}\left( C_t^{N_J} + (A^{N_J}_t + B^{N_J}_t) \varepsilon^{-3/2}\delta \right)  
=: D_t^{N_J} \sqrt{\varepsilon}.
\end{align*}
Since we assumed
\[
\delta_\varepsilon := \varepsilon^{3/2 + \eta}, \quad \text{with } 0 < \varepsilon < 1,
\]
we obtain
\begin{align}
\mathbb{E}_{J \sim \mathcal{J}}  \mathbb{E}[e_i(t)]^2
\le D_t^{N_J} \sqrt{\varepsilon}.
\end{align}
Note that $B$ is typically of order $\mathcal{O}(1)$, but could also be of order $\mathcal{O}(\varepsilon^\eta)$, which would still be compatible with the choice of $\delta$.

To obtain a high-probability bound, we apply Markov's inequality:
\[
\mathbb{P}_{J \sim \mathcal{J}} \left( 
  \mathbb{E}[e_i^2(t)] 
  \ge \tfrac{D_t^{N_J}}{\eta} \sqrt{\varepsilon}
\right)
\le \frac{
  \mathbb{E}_{J \sim \mathcal{J}} \left[ \mathbb{E}[e_i^2(t)] \right]
}{
  \frac{D_t^{N_J} \sqrt{\varepsilon}}{\eta}
}
\le \eta.
\]
\end{proof}

\subsection{Algorithms}\label{app:algorithms}
In this section, we provide pseudocodes of algorithms used in numerical experiments.
Algorithm~\ref{alg:UCB1} is the standard UCB1 algorithm stated in a general notation with the empirical mean that we estimate being the reward $\mathcal{R}_{a_{\tau}}$ defined in Equation~\ref{def:reward_function}. The algorithm is used in the sequential learning procedure to optimize $N_J$ and $\lambda$. Algorithm~\ref{alg:appr_Coverage} is a pseudocode to compute coverage $\kappa$, and Algorithm~\ref{alg:dacycle} is a pseudocode of one data assimilation cycle.

\begin{minipage}{0.48\textwidth}
\begin{algorithm}[H] 
\caption{UCB1}\label{alg:UCB1}
\begin{algorithmic}
\State \textbf{Initialization:} Play each machine once;
\For{$t=1,2, 3....$}
\State Perform action 
\[a_{t+1}=\arg \max_{a\in \mathcal{A}}  \hat{\mu}_a(t)+ \sqrt{\frac{2 \log{(t})}{N_t(a)}}\]

\State Update $ \hat{\mu}_{a+1}(t+1)$ and $N_{t+1}(a+1)$

\EndFor
\end{algorithmic}
\end{algorithm}
\end{minipage}
\begin{minipage}[h]{0.48\textwidth}
\begin{center}
\begin{tikzpicture}[decoration=brace]
    \draw(0,4)--(0,5.3);
      \draw(-0.1,4)--(0.1,4);
            \draw(-0.1,5.3)--(0.1,5.3);
            \draw[magenta!40,thick] (2,6.5) ellipse (0.4cm and 0.2cm);
            \filldraw [gray] (0,5.0) circle (2pt);
                        \filldraw [gray] (2,4.8) circle (2pt);
                          \filldraw [gray] (1,3.7) circle (2pt);
         \node at (0.3,5.0) {$\mu_1$};
               \node at (6,5.0) {\small{Optimism in the face}};
               \node at (6,4.6) {\small{of uncertainty}};
         \draw(6,5.4)--(2.2,6.5);
          \node at (1.3,3.7) {$\mu_2$};
           \node at (2.3,4.8) {$\mu_3$};
           
                    \node at (-0.3,4.65) {$\hat\mu_1$};
          \node at (1.7,5.55){$\hat\mu_3$};
           \node at (0.7,3.45) {$\hat\mu_2$};
           
                       \filldraw [magenta!40] (0,4.65) circle (2pt);
                        \filldraw [magenta!40] (2,5.55) circle (2pt);
                          \filldraw [magenta!40] (1,3.45) circle (2pt);
         
               \draw(1,2.8)--(1,4.1);
      \draw(0.9,2.8)--(1.1,2.8);
            \draw(0.9,4.1)--(1.1,4.1);

                           \draw(2,4.6)--(2,6.5);
      \draw(1.9,4.6)--(2.1,4.6);
            \draw(1.9,6.5)--(2.1,6.5);
            \node (a1) at (0, 2.3) {$a_1$};
                    \node (a2) at (1, 2.3) {$a_2$};
                            \node (a1) at (2, 2.3) {$a_3$};

  \end{tikzpicture}
  \end{center}
\end{minipage}

\begin{algorithm}[H]
  \caption{Compute approximate coverage ratio $\kappa$}\label{alg:appr_Coverage}
  \begin{algorithmic}

    \Require index set $J$, localisation radius $r_{\mathrm{loc}}$, ensemble covariance $P\in\mathbb{R}^{N_x\times N_x}$, and correlation threshold $\tau_{\mathrm{corr}}$
    \Ensure  $\kappa$  \Comment{fraction of state indices “covered”}

    \State $covered \gets 0$
    \For{$k \gets 0$ \textbf{to} $N_x-1$}
        \ForAll{$j \in J$}
            \State $d \gets \min\!\bigl(\lvert j-k\rvert,\;N_x-\lvert j-k\rvert\bigr)$
            \If{$d \le r_{\mathrm{loc}}$}
                \State $c_{kj} \gets \dfrac{\bigl(P[k,j]\bigr)^2}{\bigl(P[k,k]\bigr)^2 + 10^{-12}}$
                \If{$c_{kj} \ge \tau_{\mathrm{corr}}$}
                    \State $covered \gets covered + 1$
                    \State \textbf{break}
                \EndIf
            \EndIf
        \EndFor
    \EndFor
    \State $\kappa \gets covered / N_x$
  \end{algorithmic}
\end{algorithm}

\begin{algorithm}[H]
\caption{One Forecast-Analysis Cycle}\label{alg:dacycle}
\begin{algorithmic}
\Require $X^{\mathrm{ref}}_t, X_t, \lambda, t, \mu[\cdot], n[\cdot], c,\alpha,\beta,\gamma,\tau_{\mathrm{corr}}, r_{\mathrm{loc}},  \mathcal{J}, N_J^{\mathrm{curr}}, J^{\mathrm{curr}}$
\Ensure Updated $(X^{\mathrm{ref}}_{t+1},X_{t+1})$, diagnostics, and $(N_J^{\mathrm{next}}, J^{\mathrm{next}})$
\State Choose obs-error variance $0<\varepsilon^2<2$ \Comment{observation error order}
\State $\textit{switch} \gets \big(\mathrm{Poisson}(\lambda) > 0\big)$
\If{$N_J^{\mathrm{curr}}=\varnothing$ \textbf{or} \textit{switch}}
  \State $N_J \gets \mathrm{ChooseUCB}(\mu,n,t,c)$
  \State $J \gets$ random subset of $\{1,\dots,N_x\}$ of size $N_J$ (sorted)
\Else
  \State $N_J \gets N_J^{\mathrm{curr}}$; \quad $J \gets J^{\mathrm{curr}}$
\EndIf
\State $n[N_J] \gets n[N_J] + 1$
\State $H \gets$ row-selector for $J$; \quad $R \gets \varepsilon^2\,\mathrm{I}_{N_J}$ \Comment{define obs. operator and error covariance}
\For{$s=1$ \textbf{to} $N_{\text{Inner}}$}\Comment{$N_{\textrm{Inner}}$ steps without data assimilation}
  \State $X^{\mathrm{ref}}_{t+1} \gets \mathrm{RK4\_Step}(X^{\mathrm{ref}}_t, 0.01)$ \Comment{advance reference solution}
  \State $X \gets \mathrm{RK4\_Step}(X_t, 0.01)$ \Comment{advance ensemble, forecast step}
\EndFor
\State $\overline{X} \gets \mathrm{mean}_\mathrm{cols}(X)$; \quad $X \gets \overline{X} + 1.05(X-\overline{X})$\Comment{inflate ensemble}
\State $y \gets X^{\mathrm{ref}}_t[J] + \varepsilon\,\mathcal{N}(0,\mathrm{I}_{N_J})$ \Comment{define observations}
\State $X_{t+1} \gets \mathrm{LETKF\_GC}(y,X,H,R,r_{\mathrm{loc}})$ \Comment{update ensemble, analysis step}
\State $P \gets \frac{1}{M-1}(X_{t+1}-\overline{X}_{t+1})(X_{t+1}-\overline{X}_{t+1})^T$
\State $\kappa \gets \mathrm{Coverage}(J,r_{\text{loc}}, P, \tau_{\mathrm{corr}})$
\State $\mathcal{R}_{N_J} \gets\beta\,\kappa - \alpha\, N_J/N_x  -\gamma\,\mathrm{trace}(P)/N_x $ \Comment{reward}
\State $\mu[N_J] \gets \mu[N_J] + \dfrac{\mathcal{R}_t - \mu[N_J]}{n[N_J]}$\Comment{update the empirical mean reward}
\State $\mathrm{RMSE} \gets \|X^{\text{ref}}_{t+1}-\overline{X}_{t+1}\|$\Comment{root mean square error}
\State \Return $X^{\mathrm{ref}}_{t+1}, X_{t+1},  \mathrm{trace}(P_{t+1}), \mathrm{RMSE}, \mathcal{R}_{N_J}, \kappa, N_J, J$
\end{algorithmic}
\end{algorithm}

\subsection{Numerical experiment parameters}

\begin{table}[th]
\centering
\begin{tabular}{llp{8cm}}
\hline
\textbf{Variable} & \textbf{Setting Fig. \ref{fig:results_J_L96}} & \textbf{Description} \\
\hline
$\varepsilon$& 0.25 & Observation error order\\
$M$ & 30 & Ensemble size \\
$\Delta t$ & 0.01 & Time step for RK4 time integration \\
$N_x$ & $\{40,80,120,160\}$ & State dimension \\
$r_\textrm{loc}$ & $\{10,20,30,40\}$ & Localization radius (per $N_x$) \\
stride & $\{2,5,7,10\}$ & Spacing for allowed $N_J$ values (per $N_x$)\\
allowed\_$N_J$ & range(1,$N_x$,stride) & Candidate observation network sizes \\
$\alpha$ & 3.2 & Weight for observation count $N_J$ in reward \\
$\beta$ & 2.5 & Weight for coverage in reward \\
$\gamma$ & 0.25 & Weight for spread in reward \\
$\tau_{\text{corr}}$ & 0.30 & Threshold correlation for coverage calculation \\
ucb\_coeff & 1.0 & Coefficient in UCB exploration term \\
$\lambda$ & 1000 & Poisson intensity controlling switching frequency \\
cycles & 3000 & Forecast–analysis cycles per experiment \\
n\_reps & 50 & Number of experiment repetitions for averaging \\
$N_{\text{Inner}}$ & 5 & Number of forecast steps before DA step\\
\hline
\end{tabular}
\caption{Summary of key parameters and variables in the experiment code used for Figure \ref{fig:results_J_L96}}\label{tab:settings_runs_L96_J}
\end{table}

\end{document}